\documentclass[11pt]{amsart}
\usepackage{dsfont}
\usepackage{times}
\usepackage{enumerate}
\usepackage[euler-digits]{eulervm}
\theoremstyle{plain}
\usepackage{graphicx,color} 
\usepackage{mathtools}

\usepackage{hyperref}
\usepackage{subcaption}
\usepackage{amsmath, amssymb, graphics}
\usepackage{float}

\newtheorem{theorem}{Theorem}[section]
\newtheorem{cor}{Corollary}[section]
\newtheorem{prop}{Proposition}[section]
\newtheorem{lemma}{Lemma}[section]

\newtheorem{remark}{Remark}[section]
\newtheorem{conjecture}{Conjecture}[section]

\begin{document}

\title{Minimizing Configurations for Elastic Surface Energies with Elastic Boundaries}
\author{Bennett Palmer and \'Alvaro P\'ampano}
\date{\today}
\maketitle

\begin{abstract}
We study critical surfaces for a surface energy which contains the squared $L^2$ norm of the difference of the mean curvature $H$ and the spontaneous curvature $c_o$, coupled to the elastic energy of the boundary curve. We investigate the existence of equilibria with $H\equiv -c_o$.

When $c_o \ge 0$ we characterize those cases where the infimum of this energy is finite for topological annuli and we find the minimizer in the cases that it exists. Results for topological discs are also given.
\\

\noindent K{\tiny EY} W{\tiny ORDS}.\, Helfrich energy, bending energy, energy minimization.

\noindent MSC C{\tiny LASSIFICATION} (2010).\, 49Q10.
\end{abstract}

\section{Introduction}

Understanding the geometry of naturally occurring interfaces necessitates the consideration of an ever increasing variety of surface energy functionals accompanied by appropriate boundary energies. 

Bilipid membranes can be viewed as a fluid layer bounded by two membranes having, in general, distinct compositions and tensions. The difference in these tensions induces the layer to bend towards the interface having the lower tension \cite{LR}. This gives rise to the spontaneous curvature term in the potential energy of the bilayer which, depending on the value of a parameter $c_o$, favors the local geometry to be spherical, planar or hyperbolic. {\it We point out here that our definition for $c_o$ differs from that which is commonly used which is minus twice our quantity}.

Due to different types of stimuli, formation of pores in these membranes is quite common, resulting in an energetically unfavorable scenario (see e.g. \cite{BMF} and references therein). Therefore, bilayers usually tend to close up. Nevertheless, stabilized open bilipid membranes may exist, giving rise to equilibrium membranes with edges. These membranes have been modeled as elastic surfaces with a variety of types of boundary terms (c.f. \cite{MF}, \cite{OT1}, \cite{OT2} and \cite{WTA}). Additional previous investigations into elastic surfaces having elastic boundaries may be found in \cite{AB}, \cite{BMF} and \cite{Many}. We also call the reader's attention to the references \cite{BR}, \cite{CGS}, \cite{Tu}, \cite{Tu2}, \cite{TOY} and \cite{Z} which treat elastic surfaces with inelastic boundaries. 

In this paper we will discuss critica for compact surfaces with boundary whose principal surface energy term is the quadratic spontaneous curvature $(H+c_o)^2$, where $H$ denotes the mean curvature, and the boundary curve is treated as an elastic rod. Thus, in our case, the total energy of a surface $\Sigma$ is  \begin{equation}
\label{EHE} E[\Sigma]:=\int_\Sigma \left(a\left[H+c_o\right]^2+bK\right)d\Sigma+\oint_{\partial \Sigma}\left(\alpha\kappa^2+\beta\right)ds\:,
\end{equation}
where $a>0$, $\alpha>0$ and $c_o$, $b$, $\beta$ are any real constants. The parameter $a$ is the bending rigidity of the surface, $b$ represents its saddle-splay modulus (a name arising from the study of liquid crystals) and $c_o$ is the spontaneous curvature, while $\alpha$ is the flexural rigidity of the boundary and $\beta$ is a coupling constant. We refer to \eqref{EHE} as the \emph{Euler-Helfrich model}, continuing a naming pattern which includes the Euler-Plateau (\cite{GM}) and Kirchhoff-Plateau models (\cite{BF}, \cite{BF2}, \cite{PP}). 

Models for open bilipid membranes are often considered with an area constraint to account for fluid incompressibility. We have opted here not to include this constraint and to consider the constrained case in a future work. Obviously, any absolute minimizers for the unconstrained problem give an a priori lower bound for the constrained one. The present theoretical contribution provides a mathematical analysis of the energy $E$, paying special attention to the differential geometry of equilibria and absolute minimizers.

The Euler-Lagrange equations for $E$ consist of the fourth order differential equation
\begin{equation}\label{H}
\Delta \left(H+c_o\right)+2\left(H+c_o\right)\left(H\left[H-c_o\right]-K\right)=0\,,\nonumber
\end{equation}
together with three boundary conditions involving the curvatures and torsion of the boundary curve. If we consider the special solutions with $H+c_o\equiv 0$, then the boundary curve is necessarily an elastica circular at rest, a condition which is independent of the surface geometry. In particular, the boundary is an elastic curve when $b=0$. In the case $c_o=0$ and $b=0$, we construct some examples explicitly by numerically solving the Plateau problem.

Our main results involve the situation where the infimum of the energy within a class of surfaces occurs at the ``ground state" $H+c_o\equiv 0$. Our choice of orientation for the surface, the normal points out of a convex domain, implies that the mean curvature cannot be a positive constant and, hence, $c_o\geq 0$. The clearest results we obtain are when the surface is assumed to be a topological annulus. For example, if $b\neq 0$, and the spontaneous curvature is positive, $c_o>0$, we obtain that the infimum of the energy is finite if and only if the quantity ${\underline E}:=2\sqrt{\alpha\,\beta}-\lvert b\rvert$ is non negative, in which case the infimum is attained by a specific domain in a nodoid. When $c_o=0$ and the infimum of the energy is finite, the infimum is either attained by an axially symmetric constant mean curvature surface or it is approached asymptotically through a sequence of totally umbilical surfaces. In contrast, if $b=0$, the infimum is always finite and attained by multiple constant mean curvature surfaces. However, in this case, we cannot conclude that a minimizing surface is axially symmetric.
  
Surprisingly, the case where the surface is a topological disc is more elusive. Similar to the annular case, if $b\neq 0$ and $c_o>0$, then the infimum of the energy is finite if and only if ${\underline E}\geq 0$ holds. However, for discs there is no counterpart of the nodoidal domains, i.e. this infimum cannot be attained by constant mean curvature surfaces. Indeed, we show that for non zero spontaneous curvature, $c_o\neq 0$, and non zero saddle-splay modulus, $b\neq 0$, no constant mean curvature equilibrium can exist. Finally, when $c_o=0$ and the infimum of the energy is finite, then it is either attained or approached by totally umbilical surfaces. Similar results were obtained in \cite{BR} where the boundary was considered to be inelastic. Moreover, stability of the planar circular disc was analyzed in \cite{MF} for more general boundary energies.

\section{Euler-Helfrich Energy}

Let $\Sigma$ be a compact, connected surface with boundary and consider the immersion of $\Sigma$ in the Euclidean 3-space ${\bf R}^3$, $X:\Sigma\rightarrow {\bf R}^3\,$. Throughout this paper, we assume that $X(\Sigma)$ is an oriented surface of class $\mathcal{C}^4$ embedded in ${\bf R}^3$ with sufficiently smooth boundary $\partial\Sigma$. The unit surface normal, $\nu$, is chosen so as to point out of a convex domain. We always consider the boundary $\partial\Sigma$ as being positively oriented. 
  
For a sufficiently smooth curve $C:[0,\mathcal{L}]\rightarrow {\bf R}^3$, we denote by $s\in [0,\mathcal{L}]$ the arc length parameter of $C$. Then, if $(\,)'$ represents the derivative with respect to the arc length, the vector field $T(s):=C'(s)$ is the unit tangent to $C$. Moreover, the (Frenet) curvature of $C$, $\kappa$, is defined by $\kappa(s)=\lVert T'(s)\rVert\geq 0$. 

Each connected component of the boundary of $\Sigma$ will be represented by an arc length parameterized curve $C$. Along $C$, the oriented Darboux frame is $\{n,T,\nu\}$, where $n:=T\times \nu$ is the conormal of the boundary. The derivative of this frame with respect to the arc length parameter $s$ is given by
\begin{equation}\label{ds}
\begin{pmatrix}
n\\T\\\nu \end{pmatrix}'=
 \begin{pmatrix}
0 &-\kappa_g& \tau_g\\
\kappa_g&0&\kappa_n\\
-\tau_g&-\kappa_n&0
\end{pmatrix}\begin{pmatrix}
n\\T\\\nu \end{pmatrix},
\end{equation}
where the functions involved, $\kappa_g$, $\kappa_n$ and $\tau_g$ are, respectively, the geodesic curvature, the normal curvature and the geodesic torsion of the boundary relative to the immersion $X$. It is clear from the definitions above that $\kappa^2(s)=\lVert T'(s)\rVert^2=\kappa_g^2(s)+\kappa_n^2(s)$ holds. {\it Note that most authors define the geodesic curvature with the sign opposite to ours}. With our convention for the sign of $\kappa_g$, the classical Gauss-Bonnet Theorem reads
\begin{equation}\label{GB}
\int_\Sigma K\,d\Sigma=\oint_{\partial\Sigma}\kappa_g\,ds+2\pi\,\chi(\Sigma)\,,
\end{equation}
where $K$ is the Gaussian curvature of $\Sigma$ and $\chi(\Sigma)$ denotes its Euler-Poincar\'e characteristic.

The Euler-Helfrich functional for the immersion $X:\Sigma\rightarrow {\bf R}^3$ is the functional
\begin{eqnarray}
E[\Sigma]&:=&\int_\Sigma \left(a\left[H+c_o\right]^2+bK\right)d\Sigma+\oint_{\partial \Sigma}\left(\alpha\kappa^2+\beta\right)ds\:,\label{energy}\nonumber
\end{eqnarray}
where $a>0$, $\alpha>0$ and $c_o$, $b$, $\beta$ are any real constants. Here, $H$ denotes the mean curvature of the immersion. When necessary, we will exhibit the dependence of the functional on its parameters by writing $E= E_{a,c_o,b,\alpha,\beta}$. For convenience, we assume that all connected components of the boundary are made of the same material, so that the flexural rigidity $\alpha$ and the constant $\beta$ are the same constants for all boundary components. 

Consider sufficiently smooth variations of the immersion $X:\Sigma\rightarrow {\bf R}^3$, i.e. $X+\epsilon \delta X+{\mathcal O}(\epsilon^2)$. Below, we will denote the restriction of $\delta X$ to the boundary by $\delta C$. Then, for each term in the energy, we have the following variation formulas (for details see Appendix A):
\begin{itemize}
\item {\bf The Helfrich bending energy:}\\ 
The Helfrich energy $\mathcal{H}$, introduced in \cite{H}, is given by
\begin{equation*}
\hspace{1cm} \mathcal{H}[\Sigma]=\int_\Sigma a\left(H+c_o\right)^2\,d\Sigma\,.
\end{equation*}
The first variation of the Helfrich energy $\mathcal{H}$ for an arbitrary variation $\delta X$ is
\begin{eqnarray*}
\hspace{1cm}\delta\mathcal{H}[\Sigma]&=&a\int_\Sigma \big(\Delta H+2\left[H+c_o\right]\left(H\left[H-c_o\right]-K\right)\big)\nu\cdot \delta X\,d\Sigma\\&+&a\oint_{\partial\Sigma} \big(\left[H+c_o\right]\partial_n\left[\nu\cdot\delta X\right]-\left[\partial_n \left(H+c_o\right)\right]\nu\cdot\delta C\\&+&\left.\left[H+c_o\right]^2n\cdot\delta C\right)ds\,,
\end{eqnarray*}
where $\partial_n$ means derivative in the (outward) conormal direction along the boundary.
\item {\bf The total Gaussian curvature:}\\
The first variation of the total Gaussian curvature for $\delta X$ is given by
\begin{equation*}
\hspace{1cm}\delta\left(\int_\Sigma K\:d\Sigma\right)=\oint_{\partial \Sigma} \left(\kappa_n\partial_n\left[\nu\cdot\delta X\right]+\left[\tau_g'\nu+Kn\right]\cdot\delta C\right)ds\,.
\end{equation*}
\item {\bf The Euler bending energy:}\\  
The first variation of the bending energy of the boundary for $\delta C$ is given by
\begin{eqnarray*}
\hspace{1cm}\delta \left( \oint_{\partial \Sigma}\left[ \alpha \kappa^2+\beta\right]ds\right)=\oint_{\partial\Sigma}J'\cdot \delta C\,ds\,,
\end{eqnarray*} 
where 
\begin{equation}\label{J} 
\hspace{1cm}J:=2\alpha T''+\left(3\alpha\kappa^2-\beta\right)T\,.
\end{equation}
\end{itemize}
Combining the information above, we can express the first variation formula for the total energy $E$, by
\begin{eqnarray*}
\delta E[\Sigma]&=&\int_\Sigma a\big(\Delta H+2\left[H+c_o\right](H\left[H-c_o\right]-K)\big)\nu\cdot \delta X\,d\Sigma\\
&+&\oint_{\partial\Sigma}\left(J'+\left[a\left(H+c_o\right)^2+b K\right]n- \left[a\partial_n H-b\tau_g'\right]\nu \right)\cdot\delta C\,ds\\
&+&\oint_{\partial\Sigma}\left(a\left[H+c_o\right]+ b\kappa_n\right)\partial_n\left[\nu\cdot\delta X\right]ds\,.
\end{eqnarray*}

By considering compactly supported variations, we obtain that for an equilibrium surface
\begin{equation}\label{EL1}
\Delta H+2\left(H+c_o\right)\left(H\left[H-c_o\right]-K\right)=0
\end{equation} 
holds on $\Sigma$. Next, by taking normal variations $\delta X=\psi \nu$, we obtain the boundary integrals
\begin{equation}
0=\oint_{\partial\Sigma} \left(J'\cdot \nu -\left[a\partial_n H-b\tau_g'\right]\right)\psi\,ds+\oint_{\partial\Sigma}\left(a\left[H+c_o\right]+b\kappa_n\right)\partial_n\psi\,ds\:.\nonumber
\end{equation}
Note that $\psi$ and $\partial_n\psi$ can be chosen independently to be arbitrary functions on the boundary. By taking variations with $\psi\equiv 0$ on the boundary, we conclude that
\begin{equation}\label{EL2}
a\left(H+c_o\right)+b\kappa_n=0
\end{equation}
holds on $\partial\Sigma$. Similarly, using variations with $\partial_n\psi\equiv 0$ on $\partial\Sigma$ we also get
\begin{equation}\label{EL3}
J'\cdot \nu-a\partial_n H+b\tau_g'= 0\,,
\end{equation}
where $J$ has been defined in \eqref{J}. Finally, by taking variations that are tangential to the immersion, we also deduce the boundary condition
\begin{equation}\label{EL4}
J'\cdot n+a\left(H+c_o\right)^2+bK=0\,.
\end{equation}
Then, equations \eqref{EL1}, \eqref{EL2}, \eqref{EL3} and \eqref{EL4} are the Euler-Lagrange equations for equilibria of the total energy $E$.

Now, using the definition of $d\nu$, on the boundary $\partial\Sigma$ we have that the Gaussian curvature is given by $K:=-{\rm det}\left(d\nu\right)=\kappa_n\left(2H-\kappa_n\right)-\tau_g^2$. Hence, we plug this in equation \eqref{EL4} and use \eqref{ds} to rewrite the boundary conditions \eqref{EL3} and \eqref{EL4} respectively as
\begin{eqnarray}
\label{ELb}
2\alpha\kappa_n''+\left(\alpha\kappa^2-2\alpha\tau_g^2-\beta\right)\kappa_n+4\alpha\kappa_g'\tau_g+2\alpha\kappa_g\tau_g'+b\tau_g'-a\partial_n H&=&0\,,\nonumber \\
2\alpha\kappa_g''+\left(\alpha\kappa^2-2\alpha\tau_g^2-\beta\right)\kappa_g-4\alpha\kappa_n'\tau_g-2\alpha\kappa_n\tau_g'-b\tau_g^2&&\nonumber\\+b\kappa_n\left(2H-\kappa_n\right)+a\left(H+c_o\right)^2&=&0\,.\nonumber
\end{eqnarray}

\begin{remark} The Euler-Lagrange equations we derived above are essentially the same as equations (81)-(84) of the reference \cite{OT1}.
\end{remark}

Under a rescaling of the surface $\Sigma \rightarrow \sigma \Sigma$ for $\sigma>0$, the Willmore energy $\int_\Sigma H^2\:d\Sigma$ and the total Gaussian curvature $\int_\Sigma K\,d\Sigma$ are unchanged, the area $\int_\Sigma\,d\Sigma$ rescales quadratically and the total mean curvature $\int_\Sigma H\,d\Sigma$ and the length of the boundary rescale linearly. However, the total squared curvature of the boundary rescales like $\sigma^{-1}$, which allows for the existence of equilibrium configurations.

Using these variations by rescalings, we can prove the following result.

\begin{prop}\label{rescalings} Let $X:\Sigma\rightarrow {\bf R}^3$ be a critical immersion for the total energy $E$. Then, the following relation holds:
$$2ac_o\int_\Sigma \left(H+c_o\right)d\Sigma+\beta\mathcal{L}[\partial\Sigma]=\alpha\oint_{\partial\Sigma}\kappa^2\,ds\,,$$
where $\mathcal{L}[\partial\Sigma]$ denotes the length of the boundary $\partial\Sigma$.
\end{prop}
{\it Proof.\:} Using the information above,
\begin{eqnarray*}
E[\sigma \Sigma]&=&a\int_\Sigma H^2\:d\Sigma +b\int_\Sigma K\,d\Sigma+2ac_o\sigma\int_\Sigma H\:d\Sigma +ac_o^2\sigma^2\int_\Sigma\,d\Sigma\\
&+&\frac{\alpha}{\sigma}\oint_{\partial\Sigma}\kappa^2\,ds+\beta \sigma\mathcal{L}[\partial\Sigma]\,.
\end{eqnarray*}
Thus, differentiating with respect to $\sigma$, we get the following at the critical point ($\sigma=1$):
$$0=2ac_o\int_\Sigma H\:d\Sigma +2ac_o^2\int_\Sigma\,d\Sigma-\alpha\oint_{\partial\Sigma}\kappa^2\,ds+\beta\mathcal{L}[\partial\Sigma]\,,$$
proving the statement. {\bf q.e.d.}
\\

As a first consequence of this proposition, if either $c_o=0$ or $X:\Sigma\rightarrow{\bf R}^3$ is a critical immersion for $E$ with constant mean curvature $H=-c_o$, then
$$\beta=\frac{\alpha}{\mathcal{L}[\partial\Sigma]}\,\oint_{\partial\Sigma}\kappa^2\,ds>0$$
holds. That is, the constant $\beta$ must be positive, so from now on, in these cases, we will assume that $\beta>0$ holds.

We will see in what follows that the axially symmetric equilibria play an essential part in the study of the energy $E$ so we take some time to discuss them here. 

Let $X:\Sigma\rightarrow{\bf R}^3$ be an axially symmetric critical immersion. On parallels, which are the boundary components, both the geodesic and normal curvatures are constant, while $\tau_g\equiv 0$ holds. Moreover, in order to be a critical domain all boundary parallels must have the same radii. From this and \eqref{J} we conclude that $J=\eta T$ along the boundary parallels, where $\eta:=\alpha\kappa^2-\beta$ is a real constant.

If the surface has constant mean curvature, then by \eqref{EL3}, $\eta\kappa_n=0$ holds on $\partial\Sigma$. Assume first that $\kappa_n\equiv 0$ holds. Then combining \eqref{EL2} and \eqref{EL4} we obtain that $\eta=0$ also holds, i.e. $\eta\kappa_n=0$ if and only if $\eta=0$, which means that the radii of the boundary circles are $\sqrt{\alpha/\beta}$. 

\begin{prop}\label{axsym} Let $X:\Sigma\rightarrow{\bf R}^3$ be an axially symmetric critical immersion for $E$ and assume that on the boundary parallels the normal curvature vanishes. Then the surface has constant mean curvature $H=-c_o$ and the radii of the boundary parallels are $\sqrt{\alpha/\beta}$.
\end{prop}
{\it Proof.\:} By hypothesis, on the boundary parallels, $\kappa_g$ is a non zero constant, while $\kappa_n\equiv 0$ and $\tau_g\equiv 0$ hold. Using this in the Euler-Lagrange equations \eqref{EL2}-\eqref{EL4}, which are satisfied since the immersion is critical for $E$, we obtain
$$H+c_o=0\,,\quad\quad\quad \partial_n (H+c_o)=0\,,\quad\quad\quad \alpha\kappa^2=\beta\,,$$
on $\partial\Sigma$.

By writing the equation (\ref{EL1}) in a non parametric form and using elliptic regularity (see Theorem 6.6.1 of \cite{Mo}), one establishes the real analyticity of any $\mathcal{C}^4$ critical surface.

From (\ref{EL1}), the first two equations above and the Cauchy-Kovalevskaya Theorem, we conclude that $H+c_o\equiv 0$ holds on $\Sigma$.

Finally, the last equation above ($\alpha\kappa^2=\beta$, i.e. $\eta=0$) shows that the radii of the boundary parallels are $\sqrt{\alpha/\beta}$, as stated. {\bf q.e.d.}
\\

In the case of zero saddle-splay modulus $b$, the condition that the boundary parallels are asymptotic curves can be avoided and the result is as follows.

\begin{prop}\label{axsym0} Let $X:\Sigma\rightarrow{\bf R}^3$ be an axially symmetric critical immersion for $E$ with $b=0$. Then, either the surface has constant mean curvature $H=-c_o$ and the radii of the boundary parallels are $\sqrt{\alpha/\beta}$, or the boundary is composed by closed geodesics.
\end{prop}
{\it Proof.\:} On the boundary parallels, the geodesic and normal curvatures, $\kappa_g$ and $\kappa_n$, are constants, while $\tau_g\equiv 0$ holds. Now, since $b=0$, equation \eqref{EL2} directly gives $H+c_o=0$ on $\partial\Sigma$. Moreover, equations \eqref{EL3}-\eqref{EL4} simplify to
$$\eta\kappa_n=a\partial_n H\,,\quad\quad\quad \eta\kappa_g=0\,,$$
on $\partial\Sigma$, where $\eta:=\alpha\kappa^2-\beta$ is a constant. If $\kappa_g\equiv 0$ holds, then we are done. 

To the contrary, assume that $\kappa_g$ is a non zero constant. Then, a similar argument as in Proposition \ref{axsym} involving the Cauchy-Kovalevskaya Theorem concludes the proof.  {\bf q.e.d.}
\\

In contrast, if $b\kappa_n\neq 0$ along the boundary, then the axially symmetric critical domain $X(\Sigma)$ may not have constant mean curvature.

\section{Equilibrium Configurations with Constant Mean Curvature}

Let $C(s)$ be an arc length parameterized smooth curve in ${\bf R}^3$. We denote the Frenet frame along $C$ by $\{T,N,B\}$, where $N$ and $B$ are the unit normal and unit binormal to $C$, respectively. Note that each connected component of the boundary must be a closed curve and, consequently, their curvature can only vanish at isolated points so the Frenet frame is well defined along $C$. Moreover, the closure condition of $C$ also implies that both the curvature $\kappa(s)$ and the (Frenet) torsion $\tau(s)$ are periodic functions.

The Frenet equations involving the curvature $\kappa$ and torsion $\tau$ of a curve $C(s)$ are given by
\begin{equation}\label{feq}
\begin{pmatrix}
T\\N\\B \end{pmatrix}'=
 \begin{pmatrix}
0 &\kappa& 0\\
-\kappa&0&\tau\\
0&-\tau&0
\end{pmatrix}\begin{pmatrix}
T\\N\\B \end{pmatrix},
\end{equation}
where $\left(\,\right)'$ means derivative with respect to the arc length parameter $s$. 

Denote by $\theta$ the oriented angle between the normal to the boundary component $C$, $N$, and the normal to the surface $\Sigma$, $\nu$. This angle $\theta\in\left[-\pi,\pi\right]$ will be referred to as the contact angle between the surface and the boundary. (This angle is actually the contact angle between the surface and the ruled, developable  surface given by $(s,t)\mapsto T(s)+tN(s)$.) Write
\begin{equation}\label{cframe}
\nu+i n=e^{i\theta}\left(N+i B\right),
\end{equation}
then using \eqref{ds} and \eqref{feq} we get
\begin{eqnarray}
\kappa_g&=&\kappa\sin\theta\,,\label{kg}\\
\kappa_n&=&\kappa\cos\theta\,, \label{kn}\\
\tau_g&=&\theta'-\tau\,.\label{taug}
\end{eqnarray}

Consider a critical immersion with constant mean curvature $H$. If $H\neq -c_o$, it is clear from the Euler-Lagrange equation \eqref{EL1} that the Gaussian curvature must also be constant and satisfy 
\begin{equation}\label{Ko}
K\equiv H\left(H-c_o\right).
\end{equation}
Hence, the surface is isoparametric, so it must  be a compact domain of a plane, a sphere or a right circular cylinder. We then have the following result.

\begin{theorem}\label{clas} Let $X:\Sigma \rightarrow {\bf R}^3$ be a critical immersion for $E$ with constant mean curvature $H\neq -c_o$. Then, the surface is a compact domain in the sphere of radius $1/\lvert H\rvert$ bounded by disjoint circles of radii $\sqrt{\alpha/\beta}$. Moreover, the energy parameters must verify $a=-b>0$, $c_o=0$ and $H^2\leq \beta/\alpha$.
\end{theorem}
{\it Proof.\:} As mentioned above, if $H\neq -c_o$ equation \eqref{EL1} implies that $X(\Sigma)$ must be a domain in an isoparametric surface. We study each case separately.

If $X(\Sigma)$ is a planar domain, then $H=K=0$ and we get $\kappa_n=\tau_g=0$ along the boundary. Plugging this in the Euler-Lagrange equation \eqref{EL2} we obtain that $ac_o=0$, but this is not possible, since we are assuming $a> 0$ and $0=H\neq -c_o$ so there are no planar domains critical for $E$ with $H\neq -c_o$.

Next, suppose that $X(\Sigma)$ is a compact domain in a right circular cylinder. Cylinders are flat and have constant non zero mean curvature, thus it follows from equation \eqref{Ko} that $H=c_o\neq 0$. Using this in \eqref{EL2}, we get that $2ac_o+b\kappa_n=0$. If $b=0$ we obtain once more $ac_o=0$, which is not possible. Thus, $b\neq 0$ and $\kappa_n\equiv{\rm constant}\neq 0$. Moreover, along the boundary, we also have
$$0=K=\kappa_n\left(2H-\kappa_n\right)-\tau_g^2\,,$$
which implies that $\tau_g$ is constant too. Then, we can simplify the Euler-Lagrange equations \eqref{EL3} and \eqref{EL4}, obtaining
\begin{eqnarray}
\left(\alpha\kappa_g^2+\eta\right)\kappa_n+4\alpha\kappa_g'\tau_g&=&0\,,\label{dem1}\\
2\alpha\kappa_g''+\left(\alpha\kappa_g^2+\eta\right)\kappa_g+4ac_o^2&=&0\,, \label{dem2}
\end{eqnarray}
where $\eta:=\alpha\kappa_n^2-2\alpha\tau_g^2-\beta\equiv{\rm constant}$. 

If $\tau_g\equiv 0$, equation \eqref{dem1} reduces to $\kappa_g^2=-\eta/\alpha\equiv{\rm constant}$. This can be used in \eqref{dem2} to get $4ac_o^2=0$, which is a contradiction. Next, we suppose $\tau_g\neq 0$. Differentiating \eqref{dem1} and substituting in \eqref{dem2} we conclude once more that $\kappa_g$ must be constant, but this is impossible since again $4ac_o^2\neq 0$. That is, there are no critical domains for $E$ with $H\neq -c_o$ in right circular cylinders.

Finally, we study the case where $X(\Sigma)$ is a spherical domain. The sphere is totally umbilical, i.e. $K=H^2\neq 0$. From \eqref{Ko} we get $c_o=0$. Then, as before, from \eqref{EL2}, we get $aH+b\kappa_n=0$. Again, if $b=0$ necessarily $aH=0$, which is a contradiction. Therefore, $b\neq 0$ and $\kappa_n=-aH/b\neq 0$. Moreover, due to the definition of the Gaussian curvature along $\partial\Sigma$ we also get that $\tau_g$ is constant and
$$\tau_g^2=\kappa_n\left(2H-\kappa_n\right)-H^2=-\frac{\left(a+b\right)^2}{b^2} H^2\,.$$
Since $H\neq 0$, this is only possible if $\tau_g=0$ and $a=-b$. Note that equations \eqref{dem1} and \eqref{dem2} remain true (with $c_o=0$). For $\tau_g=0$, \eqref{dem1} and \eqref{dem2} simplify to $\alpha\kappa^2=\beta$. Thus, $\kappa$ must be constant. This together with $\kappa_n$ being constant implies that $\theta$ is also constant, \eqref{kn}. Next, we use \eqref{taug} to conclude that $0=\tau_g=-\tau$, i.e. we have a domain in the sphere of radius $1/\lvert H\rvert$ bounded by circles of the same radii, namely, $\sqrt{\alpha/\beta}$. 

Note that in order for the boundary circles to lie on the sphere, their radii must be smaller or equal the radius of the sphere, i.e. we conclude that $H^2\leq \beta/\alpha$. This finishes the proof. {\bf q.e.d.}
\\

Since the energy coefficients satisfy $a=-b>0$ and $c_o=0$ in the case given in the previous theorem, we have that spherical domains are, precisely, critical immersions for the Willmore energy with elastic boundary, i.e. ($W=E_{a,c_o=0,b=-a,\alpha,\beta}$)
\begin{equation}\label{Willmoreelastic}
W[\Sigma]:=a\int_{\Sigma}\left(H^2-K\right)d\Sigma+\oint_{\partial\Sigma}\left(\alpha\kappa^2+\beta\right)ds\,,
\end{equation}
with $a$, $\alpha$ and $\beta$ positive constants. As a consequence of Proposition \ref{rescalings}, since $c_o=0$, then $\beta>0$ holds.

For non zero, sufficiently small values of $H$ we can remove $l$ disjoint discs $D_i$ of radii $\sqrt{\alpha/\beta}$ from the sphere of radius $1/\lvert H\rvert$ and obtain a critical point of $W$ having energy
$$W\left[S^2\setminus \bigcup_{i=1}^{l}D_i\right]=\oint_{\partial\Sigma}\left(\alpha\kappa^2+\beta\right)ds=4\pi l\sqrt{\alpha\,\beta}\,.$$

If $X:\Sigma\rightarrow{\bf R}^3$ is an immersion critical for the energy $E$ with constant mean curvature $H=-c_o$, then \eqref{EL1} trivially holds. In this case, as a consequence of Proposition \ref{rescalings} the coefficient $\beta$ is positive. Assuming first that $b\neq 0$, then \eqref{EL2}, \eqref{EL3} and \eqref{EL4} simplify, respectively, to
\begin{eqnarray}
\kappa_n&=&0\,,\label{el2}\\
4\alpha\kappa_g'\tau_g+2\alpha\kappa_g\tau_g'+b\tau_g'&=&0\,,\label{el3}\\
2\alpha\kappa_g''+\left(\alpha\kappa_g^2-2\alpha\tau_g^2-\beta\right)\kappa_g-b\tau_g^2&=&0\,.\label{el4}
\end{eqnarray}
Note that, since $\kappa\neq 0$, by combining \eqref{el2} and \eqref{kn}, we conclude that the contact angle $\theta$ is constant, $\theta\equiv \pm\pi/2$. This implies, using \eqref{kg} and \eqref{taug}, that $\tau_g=-\tau$ and $\kappa_g=\pm\kappa$. Then, equations \eqref{el3} and \eqref{el4} can be expressed in terms of the curvature and torsion as
\begin{eqnarray}
4\alpha\kappa'\tau+2\alpha\kappa\tau'\pm b\tau'&=&0\,,\label{elc1}\\
2\alpha\kappa''+\left(\alpha\kappa^2-2\alpha\tau^2-\beta\right)\kappa\mp b\tau^2&=&0\,.\label{elc2}
\end{eqnarray}

If $b=0$, then not only \eqref{EL1} but also equation \eqref{EL2} holds. Moreover, from the Euler-Lagrange equations \eqref{EL3} and \eqref{EL4}, we get that $J'\cdot n=J'\cdot \nu=0$ where $J$ is given in \eqref{J}. This means $J'$ has no component normal to the curve.

Using \eqref{J} and \eqref{feq}, we find $J=(\alpha \kappa^2-\beta)T+2\alpha \kappa'N+2\alpha \kappa \tau B$. Differentiating this using \eqref{feq} then gives $J'\cdot T=0$. We conclude that $J'\equiv 0$ along the boundary $\partial\Sigma$. Now, rewriting $J'$ in terms of the Frenet frame, we see $J'\equiv 0$ is equivalent to equations \eqref{elc1} and \eqref{elc2} when $b=0$. 

Equations \eqref{elc1} and \eqref{elc2} for any $b\in{\bf R}$ are the Euler-Lagrange equations of a curvature energy, as we show in the following proposition.

\begin{prop}\label{criticality} Let $X:\Sigma\rightarrow {\bf R}^3$ be an equilibrium immersion for $E$ with constant mean curvature $H=-c_o$, then each connected component of $\partial\Sigma$ is a simple closed critical curve for
\begin{equation}\label{energyc}
F[C]\equiv F_{\mu,\lambda}[C]:=\int_C\left( \left[\kappa+\mu\right]^2+\lambda\right)ds
\end{equation}
in ${\bf R}^3$, where the energy parameters are given by 
$$\mu:=\pm\frac{b}{2\alpha}\,,\quad\quad\quad \lambda:=\frac{\beta}{\alpha}-\mu^2\:.$$
\end{prop}
{\it Proof.\:} Consider the energy $F$ acting on the space of closed curves. By standard arguments involving integration by parts (for details adapt the computations of Appendix A), we have
$$\delta F[C]=\oint_{C} \bar{J}'\cdot \delta C\,ds$$
where the vector field $\bar{J}$ is defined along $C$ by
\begin{equation}\label{Jbar}
\bar{J}:=\left(\kappa^2-\left[\lambda+\mu^2\right]\right)T+2\kappa'N+2\tau\left(\kappa+\mu\right) B\, ,
\end{equation}
so $\bar{J}'=0$ is the Euler-Lagrange equation for $F$.

If we consider now the values of the parameters $\mu$ and $\lambda$ as in the statement, it turns out that the tangent component of $\bar{J}'=0$ is an identity, while its binormal and normal components are \eqref{elc1} and \eqref{elc2}, respectively. Therefore, any closed connected component of the boundary satisfies the Euler-Lagrange equations and, hence, it is a critical curve for $F$, proving the result. {\bf q.e.d.}
\\

The energy \eqref{energyc} is an extension of the classical bending energy of curves, which appears when $\mu=0$ (i.e. $b=0$). It is referred to as the \emph{bending energy circular at rest} since, if $\mu\neq 0$, it can be used to study the shape of stiff rods which are circular in their undeformed state, \cite{CCG}. Of course, the parameter $\mu$ is just an analogue for curves of what the spontaneous curvature is for surfaces. 

Note that, since the total curvature is constant on a regular homotopy class of planar curves, the only planar critical curves for $F$ are Euler's elasticae, \cite{E}. In particular, among closed ones, only circles and elastic figure eights appear, the former being the only embedded ones.  

Modifying the computations of \cite{LS} which were carried out for Kirchhoff rods, we obtain, along any critical curve for $F$,  two commuting Killing vector fields given by $\bar{J}$, \eqref{Jbar}, and 
\begin{equation}\label{Ibar}
\bar{I}:=2\left(\kappa+\mu\right)B\,.
\end{equation}
These vector fields can be used to integrate the Euler-Lagrange equations for $F$ once. It can be shown that the length of the vector field $\bar{J}$ is constant along critical curves, in particular $\lVert \bar{J}\rVert^2=d\geq 0$. Moreover, the inner product of $\bar{J}$ and $\bar{I}$ is also constant with  $\bar{J}\cdot\bar{I}=e\in{\bf R}$. Hence, we obtain the first integrals
\begin{eqnarray}
4\tau\left(\kappa+\mu\right)^2&=&e\,,\label{fi1}\\
4\left[\kappa'\right]^2+\left(\kappa^2-\left[\lambda+\mu^2\right]\right)^2+4\tau^2\left(\kappa+\mu\right)^2&=&d\,,\label{fi2}
\end{eqnarray}
where $d\geq 0$ and $e$ are the constants of integration.

If $d=0$, from \eqref{fi2} we get that $\kappa^2=\lambda+\mu^2$ is constant and that $\tau\left(\kappa+\mu\right)=0$. Plugging the last equation in \eqref{fi1} we conclude that $e=0$. In particular, if $\kappa\neq -\mu$ is constant, then $\tau=0$ holds and the critical curve $C$ is a circle of radius $\sqrt{\alpha/\beta}$.

Now consider $d>0$. Following \cite{LS}, the vector fields $\bar{J}$ and $\bar{I}$ can be uniquely extended to Killing vector fields defined on the whole space ${\bf R}^3$. For convenience, we also denote these extensions by $\bar{J}$ and $\bar{I}$, respectively. Let $r$, $\vartheta$, $z$ be cylindrical coordinates. Since $\bar{J}$ has constant length, it is a translational vector field. After rigid motions, we may assume that $\bar{J}=\sqrt{d}\,\partial_z$. With respect to the same coordinates, $\bar{I}$ can be expressed as (see \cite{LS} for details)
\begin{equation}\label{Iexp}
\bar{I}=\sqrt{d}\,\partial_\vartheta+\frac{e}{\sqrt{d}}\partial_z\,.
\end{equation}
Since $\lVert\partial_\vartheta\rVert^2=r^2(s)$, using equations \eqref{Ibar} and \eqref{Iexp} to compute $\lVert\bar{I}\rVert^2=\bar{I}\cdot\bar{I}$, we conclude that
\begin{equation}\label{r(s)}
r^2(s)=\frac{1}{d^2}\left(4d\left[\kappa+\mu\right]^2-e^2\right).
\end{equation}
For a critical curve $C$, we write its unit tangent vector field as $T(s)=r'(s)\partial_r+\vartheta'(s)\partial_\vartheta+z'(s)\partial_z$. Hence, from \eqref{Jbar} and $\bar{J}=\sqrt{d}\,\partial_z$, we compute $T\cdot \bar{J}$, obtaining that
\begin{equation}\label{zprima}
z'(s)=\frac{1}{\sqrt{d}}\left(\kappa^2-\left[\lambda+\mu^2\right]\right).
\end{equation}
Finally, we use equation \eqref{Ibar} to get $0=T\cdot \bar{I}$, which gives rise to
\begin{equation}\label{thetaprima}
\vartheta'(s)=-e\sqrt{d}\,\frac{\kappa^2-\left(\lambda+\mu^2\right)}{4d\left(\kappa+\mu\right)^2-e^2}=-\frac{e\,z'(s)}{d\,r^2(s)}\,.
\end{equation}
In conclusion, if $C$ is a critical curve for $F$ with $d>0$, it can be parameterized using cylindrical coordinates as $C(s)=\left(r(s)e^{i\vartheta(s)},z(s)\right)$ where $r(s)$, $z(s)$ and $\vartheta(s)$ are obtained from \eqref{r(s)}, \eqref{zprima} and \eqref{thetaprima}, respectively. All of the functions involved in the parameterization of $C$ are described in terms of its curvature, $\kappa(s)$. By the Fundamental Theorem of Curves, the curvature $\kappa(s)$ and the torsion $\tau(s)$ completely determine, up to rigid motions, an arc length parameterized curve $C(s)$. Due to equation \eqref{fi1}, for $d>0$, the torsion of critical curves can be expressed in terms of the curvature $\kappa(s)$. Hence, up to rigid motions, critical curves with $d>0$ are completely described by their curvature $\kappa(s)$, which is a solution of \eqref{fi2}.

The periodicity of the curvature and torsion of the boundary components is a necessary but not sufficient condition for  the curve to be closed. Let $\varrho$ be the period of the curvature of a critical curve $C$. From equation \eqref{fi1}, the torsion of $C$ is also periodic of period $\varrho$. Then, from the parameterization given above, the closure conditions are given by the following:
\begin{eqnarray*}
\Delta z&=&\sqrt{d}\int_0^\varrho z'(s)\,ds=\int_0^\varrho \left(\kappa^2-\left[\lambda+\mu^2\right]\right)ds=0\,,\\ 
\Delta\vartheta&=&\int_0^\varrho \vartheta'(s)\,ds=-e\sqrt{d}\int_0^\varrho \frac{\kappa^2-\left(\lambda+\mu^2\right)}{4d\left(\kappa+\mu\right)^2-e^2}\,ds=2\frac{p}{q}\pi\,,
\end{eqnarray*}
for natural numbers $p$ and $q$ satisfying that ${\rm gcd}(p,q)=1$. These parameters have a geometric meaning. The integer $p$ represents the number of times the critical curve winds around the $z$-axis while $q$ is the number of periods of the curvature contained in one period of the curve.

Given periodic functions $\kappa$ and $\tau$ that are solutions of \eqref{fi1}-\eqref{fi2}, the closure conditions above are difficult to check analytically. In the special case $b=0$ (i.e. $\mu=0$), the curvature can be explicitly obtained with the aid of elliptic integrals and these conditions were studied in \cite{ls}, proving that there exist infinitely many embedded closed non planar elastic curves lying on rotational tori. These non planar elastic curves represent $(q,p)$-torus knots for $2p<q$, denoted here by $G(q,p)$. From formulas \eqref{r(s)} and \eqref{zprima}, it is clear that for an arbitrary $\mu$, simple closed critical curves with non constant curvature also lie on rotational tori, representing $(q,p)$-torus knots $G(q,p)$. 

From Proposition \ref{criticality} and the argument above, the boundary of a critical immersion for $E$ with $H\equiv -c_o$ is composed of critical curves for $F$ either having constant curvature or representing a torus knot. As a consequence, in the latter case, $\Sigma$ is a Seifert surface for the link $\partial\Sigma$ (see e.g. \cite{M}).

Using properties of Seifert surfaces, we obtain the following topological restriction for the existence of critical immersions for $E$ with constant mean curvature $H=-c_o$.

\begin{prop}\label{genus} Let $X:\Sigma\rightarrow{\bf R}^3$ be an immersion of a genus $g$ surface $\Sigma$ with constant mean curvature $H=-c_o$ critical for $E$. Assume that $\Sigma$ has only one boundary component $C$. Then, either $C$ has constant curvature or 
$$2g\geq \left(p-1\right)\left(q-1\right)$$
holds, where $p$ and $q$ are the natural numbers describing the type of the torus knot $C\cong G(q,p)$. 
\end{prop}
{\it Proof.\:} Since the critical immersion $X:\Sigma\rightarrow{\bf R}^3$ has constant mean curvature $H=-c_o$, its boundary component is a simple closed critical curve for $F$. If it has non constant curvature, then it represents a $(q,p)$-torus knot, $G(q,p)$, for $p,q\geq 1$, and $\Sigma$ is a Seifert surface.

The genus of a knot $G$, $g(G)$, is defined as the minimum of the genus of any of its Seifert surfaces. In particular, for a $(q,p)$-torus knot, $G\equiv G(q,p)$, with $p,q\geq 1$ its genus can be explicitly computed (see \cite{M} for details). Thus, we have
$$\frac{1}{2}\left(p-1\right)\left(q-1\right)=g(G)\leq g$$
where $g$ is the genus of $\Sigma$. {\bf q.e.d.}
\\

In the case $b=0$, we can construct equilibrium surfaces with $H\equiv -c_o\leq 0$ by starting with a simple closed elastic curve $C$ in ${\bf R}^3$. Hildebrandt's solution of the volume constrained Plateau Problem \cite{Hil} can then be applied to obtain, for $c_o$ sufficiently small, a topological disc with constant mean curvature $H=-c_o$ having boundary curve $C$ (see also \cite{Heinz}). Note that these surfaces may have branch points. Indeed, it follows from Proposition \ref{genus} that this is the case if $C$ is a knotted elastic curve.

\begin{figure}[h!]
\centering
\begin{subfigure}[b]{0.295\linewidth}
\includegraphics[width=\linewidth]{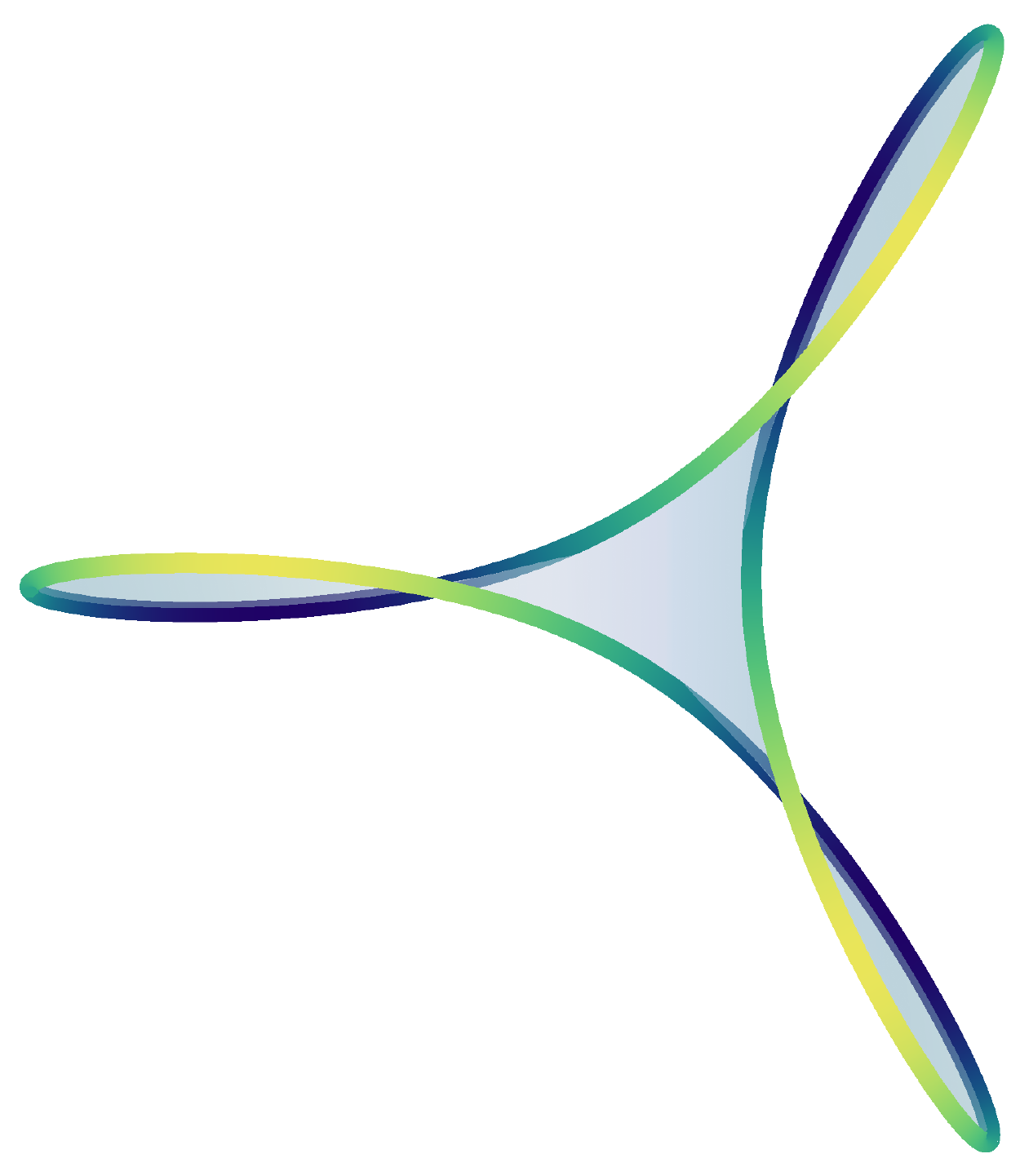}
\end{subfigure}
\begin{subfigure}[b]{0.34\linewidth}
\includegraphics[width=\linewidth]{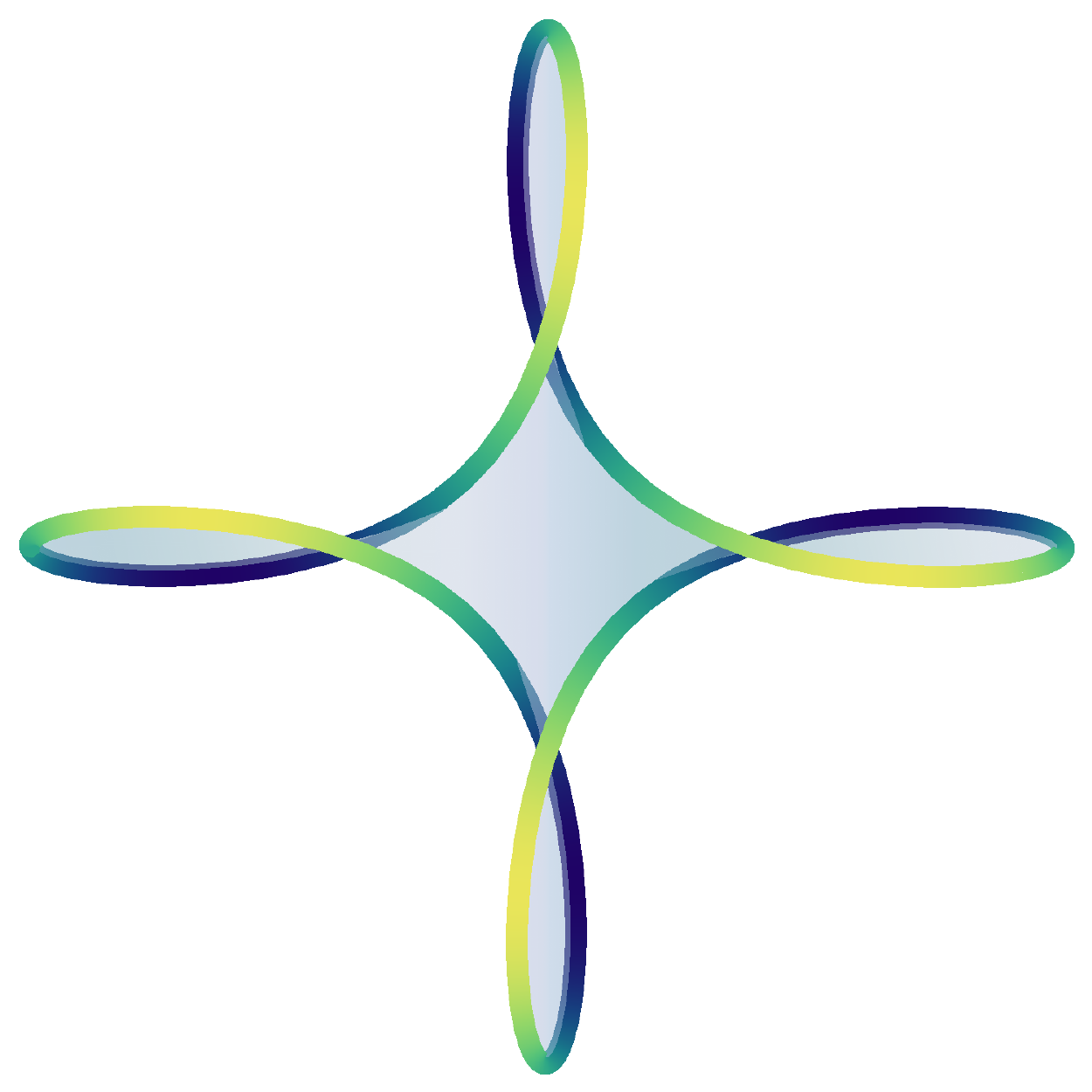}
\end{subfigure}
\,
\begin{subfigure}[b]{0.322\linewidth}
\includegraphics[width=\linewidth]{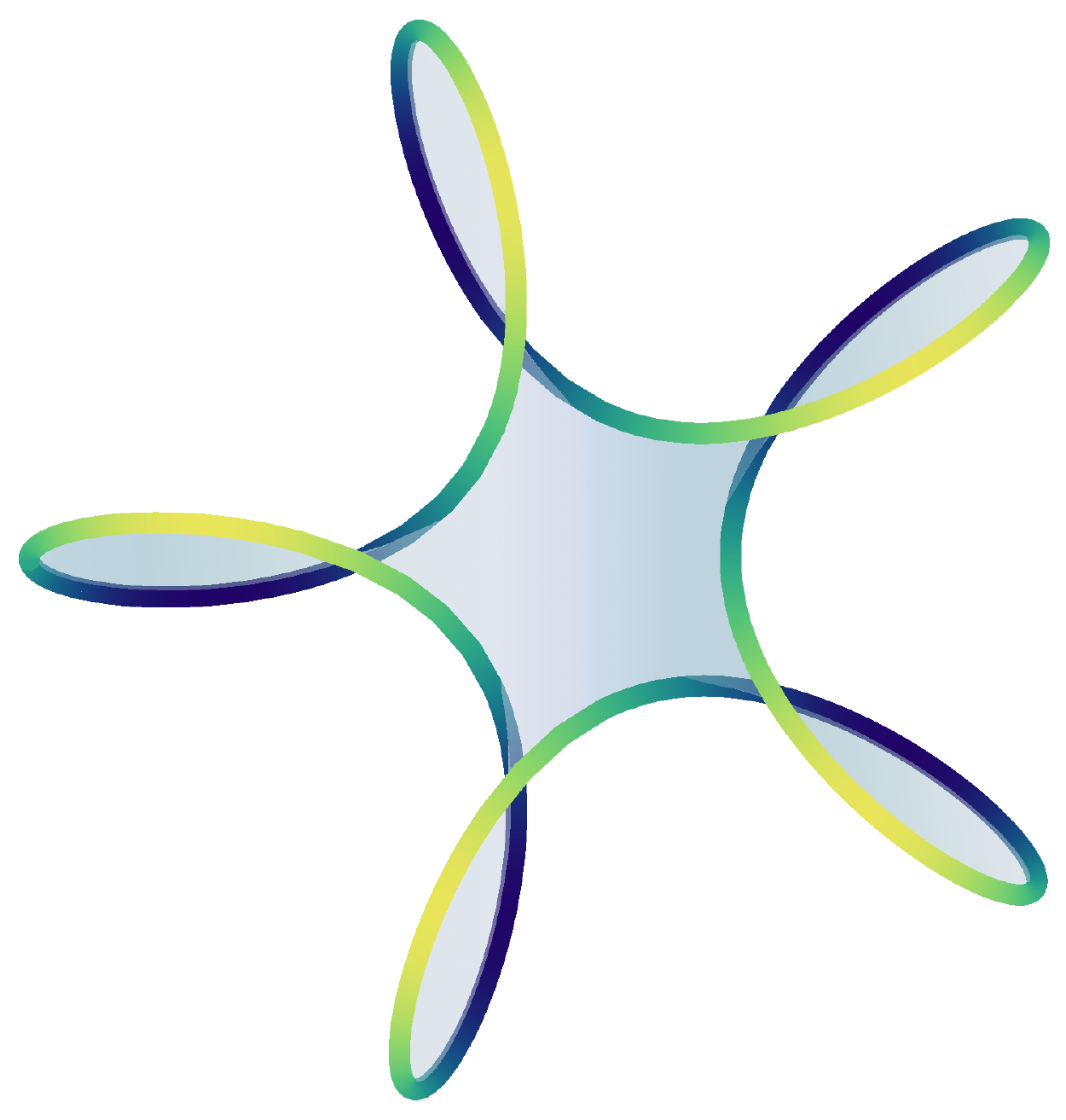}
\end{subfigure}
\caption{Minimal surfaces of disc type spanned by elastic curves of type $G(q,1)$ for $q=3, 4, 5$. These configurations are critical for $E$ with $b=0$ and $c_o=0$, i.e. for $E_{a,c_o=0,b=0,\alpha,\beta}$.}\label{minelas}
\end{figure}

In Figure \ref{minelas}, we show some minimal embedded disc type surfaces whose boundaries are elastic curves of type $G(q,1)$. These minimal discs have been obtained numerically by using \emph{Wolfram Mathematica} to implement a numerical version of the mean curvature flow (\cite{G}) to solve the Plateau Problem for a fixed elastic boundary curve.

\begin{figure}[h!]
\centering
\begin{subfigure}[b]{0.325\linewidth}
\includegraphics[width=\linewidth]{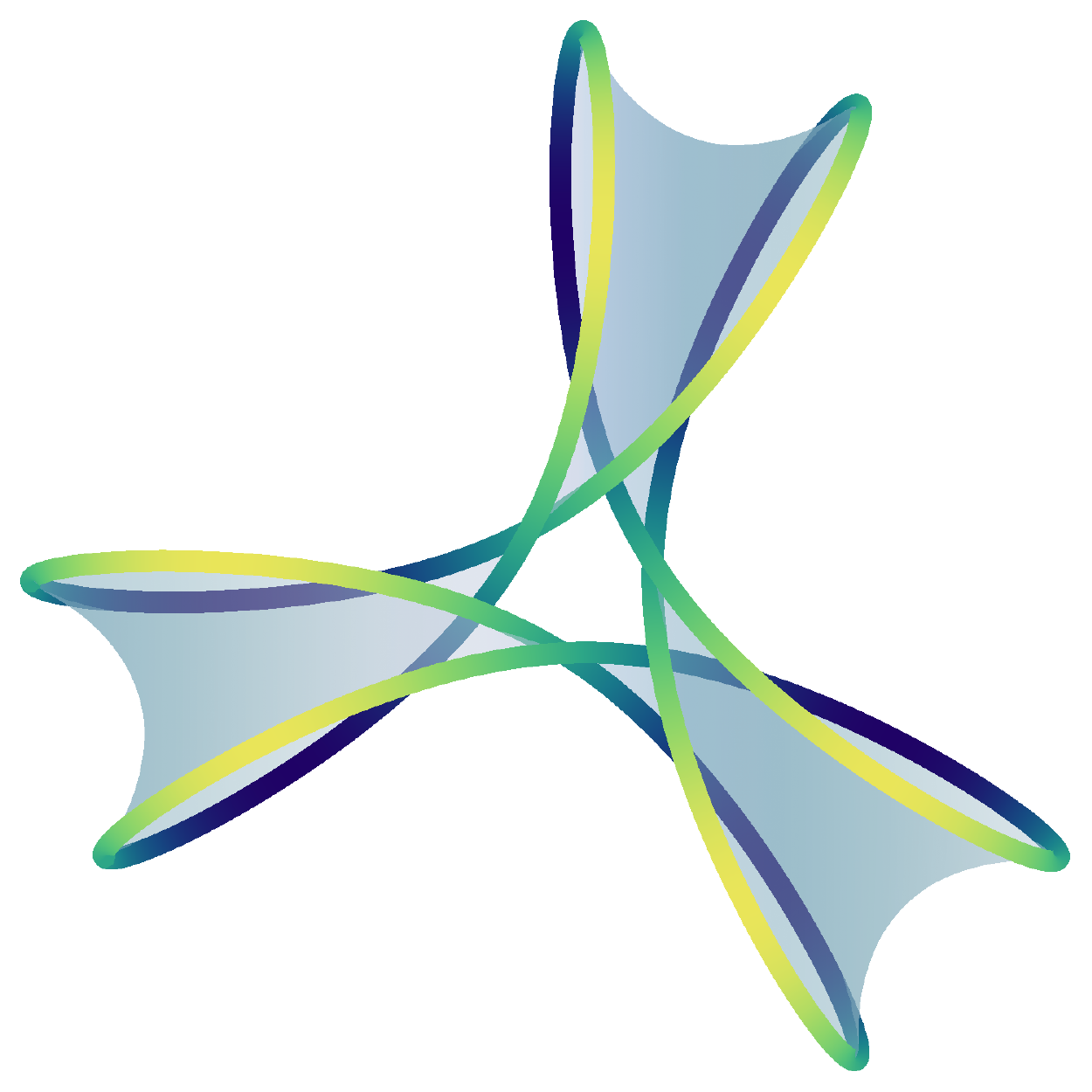}
\caption{$G(3,1)$}
\end{subfigure}
\begin{subfigure}[b]{0.325\linewidth}
\includegraphics[width=\linewidth]{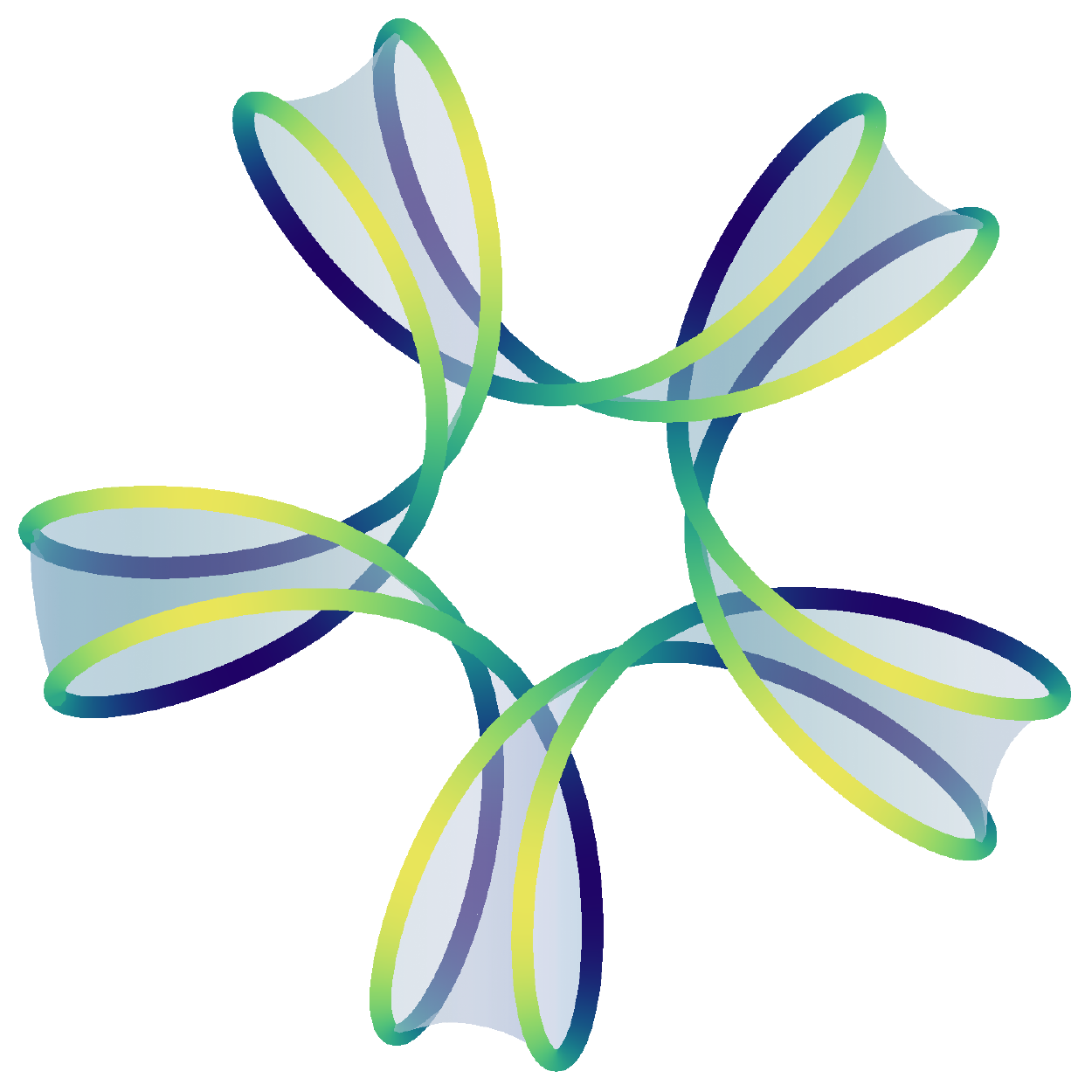}
\caption{$G(5,1)$}
\end{subfigure}
\begin{subfigure}[b]{0.325\linewidth}
\includegraphics[width=\linewidth]{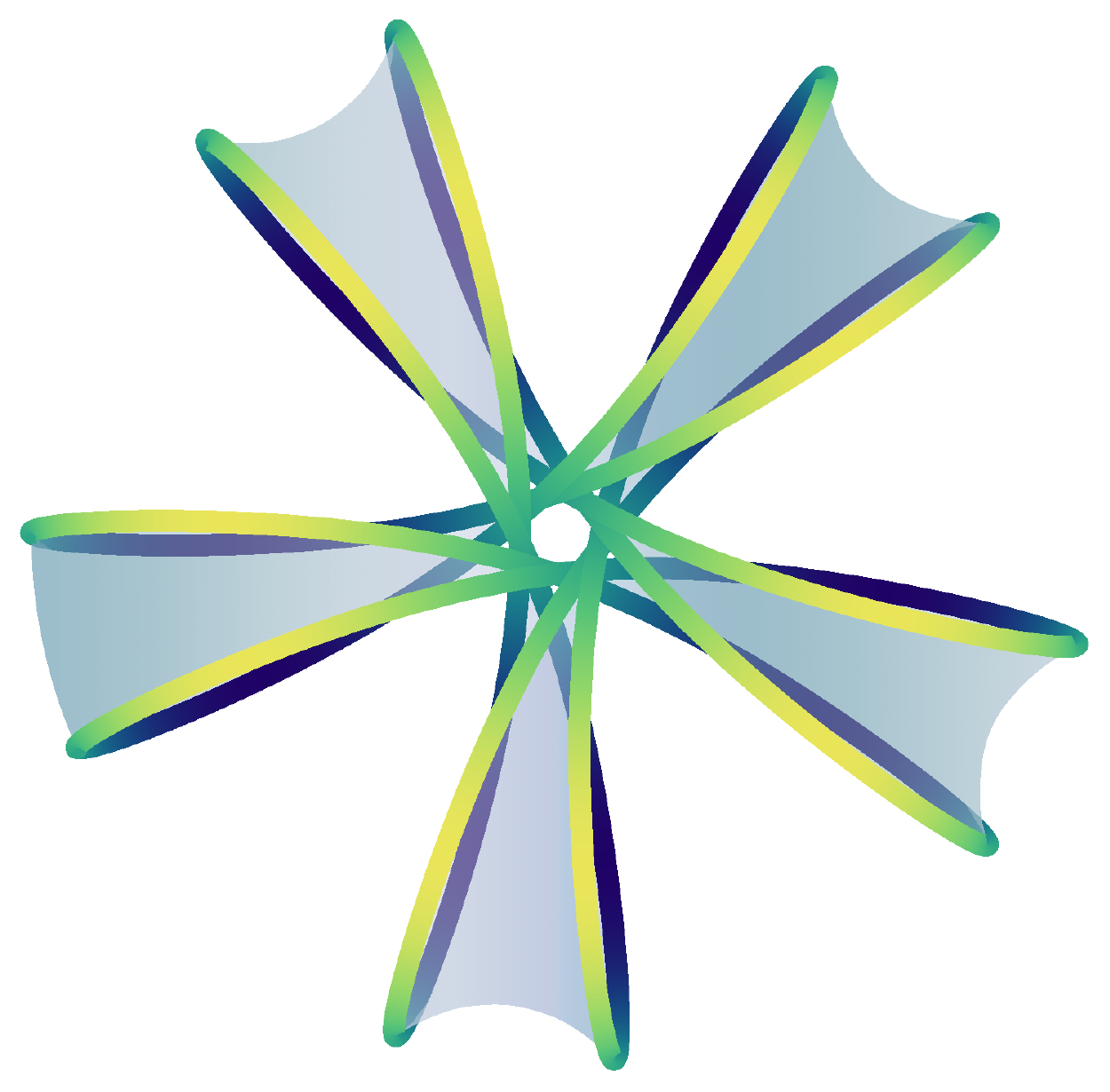}
\caption{$G(5,2)$}
\end{subfigure}
\caption{Minimal annuli with boundary composed by two congruent elastic curves representing torus knots $G(q,p)$. These domains are critical for $E_{a,c_o=0,b=0,\alpha,\beta}$.}\label{minelasan}
\end{figure}

The same numerical algorithm can be applied to produce minimal, critical annuli when $b=0$. In Figure \ref{minelasan}, we show some annuli obtained in this way. We point out that the restriction of Proposition \ref{genus} that the surface has only one boundary component is essential. For instance, in Figure \ref{minelasan} (C), a minimal annulus bounded by two congruent elastic curves of type $G(5,2)$ is shown.

Producing examples in the case $b\neq 0$ is substantially more involved. In this case, the surface is more strongly coupled to the boundary curve by the condition $\kappa_n\equiv 0$ on $\partial\Sigma$, i.e. each boundary component is a closed asymptotic curve. Although globally, the existence of surfaces with $H\equiv -c_o$ genus zero critical for $E$ with suitable prescribed boundary is problematic, locally, it is always possible to construct such surfaces. Let $\kappa(s)$ and $\tau(s)$ be solutions of \eqref{fi1}-\eqref{fi2} defined on a sufficiently small interval $I$. Since the coefficients of \eqref{fi1} and \eqref{fi2} are real analytic, then both $\kappa(s)$ and $\tau(s)$ are real analytic functions of the arc length $s$. By the Fundamental Theorem of Curves, there exists a curve $C(s)$ defined on $I$, which can be found by solving the Frenet equations \eqref{feq}. By the same argument as was used above, $C(s)$ is also real analytic on $I$. 

Define a unit vector field $\nu(s)$ along $C(s)$ (in the case $b\neq 0$, we ask that $\nu(s)$ is the unique unit vector orthogonal to $T(s)$ making an angle $\theta\equiv\pm\pi/2$ with the Frenet normal $N(s)$, that is, from \eqref{cframe}, $\nu(s)\equiv\pm B(s)$). By analyticity, both the curve $C(s)$ and the vector field $\nu(s)$ have holomorphic extensions $C(z)$ and $\nu(z)$, for a complex variable $z=s+i t$, to a simply connected domain $U$ in the complex plane ${\bf C}$, with $I\subset U$.

Next, for fixed $s_o\in I$, the Bj\"{o}rling's Formula
$$X(z):=\Re\left(C(z)+i\int_{s_o}^z \left[C'(\omega)\times \nu(\omega)\right]d\omega\right)$$
gives a minimal surface containing the curve $C$ having unit normal $\nu(s)$ along $C$. Finally, if we take $\Sigma$ to be the part of the minimal surface $X(z)$ which lies on one side of the curve $C$, we obtain a ``local'' equilibrium configuration for $E$. Note that from our definition of $\nu(s)$ for the case $b\neq 0$, $\kappa_n\equiv 0$ holds along $C$.

This method allows us to locally construct minimal surfaces. For non zero constant mean curvature $H=-c_o\neq 0$, locally, it is also possible to construct these surfaces via the loop group formulation, \cite{BD}. 

To conclude this section, we are going to prove that the condition $H\equiv -c_o$ on $\Sigma$ can be deduced as a consequence of the surface making a suitable constant contact angle $\theta$ along a planar boundary component.

\begin{prop}\label{cca} Let $X:\Sigma\rightarrow{\bf R}^3$ be a critical immersion for $E$. If any connected component of $\partial\Sigma$ is planar and its contact angle satisfies $\theta\equiv \pm\pi/2$, then the surface has constant mean curvature $H=-c_o$.
\end{prop}
{\it Proof.\:} Let $C$ be a planar connected component of $\partial\Sigma$ and assume that $\theta\equiv\pm\pi/2$ holds along $C$. Then, using equations \eqref{kg}-\eqref{taug}, we get that $\kappa_g=\pm\kappa$, $\kappa_n\equiv 0$ and $\tau_g\equiv 0$ hold along $C$. Therefore, the Euler-Lagrange equation \eqref{EL2} reduces to $H=-c_o$ on $C$.

Moreover, substitution of these quantities in the Euler-Lagrange equation \eqref{EL3} implies that
$$J'\cdot \nu+b\tau_g'=0=a\partial_n H$$
holds on $C$. This together with \eqref{J}, \eqref{EL3} and \eqref {EL4}, implies that $C$ is a real analytic curve.

Since the surface is a solution of the Euler-Lagrange equation \eqref{EL1} with $H+c_o=0$ and $\partial_n(H+c_o)=0$ on a boundary component, applying the Cauchy-Kovalevskaya Theorem as in Section 2, we conclude that $X(\Sigma)$ has constant mean curvature $H=-c_o$. {\bf q.e.d.}
\\

Observe that Proposition \ref{axsym} can be deduced from previous result. Indeed, for axially symmetric immersions the boundary components are planar and from \eqref{kn}, $\kappa_n\equiv 0$ is equivalent to $\theta\equiv\pm\pi/2$. 

In the case where the saddle-splay modulus $b$ vanishes, Proposition \ref{cca} can be improved as follows.

\begin{prop}\label{ccab0} Let $X:\Sigma\rightarrow{\bf R}^3$ be a critical immersion for $E$ with $b=0$. If any connected component of $\partial\Sigma$ is planar and its contact angle is a constant satisfying $\theta\neq 0$, $\pm\pi$, then the surface has constant mean curvature $H=-c_o$.
\end{prop}
{\it Proof.\:} As before, we denote by $C$ the planar connected component of $\partial\Sigma$ and we assume that $\theta\neq 0$, $\pm\pi$ is constant along $C$. From \eqref{kg}-\eqref{taug}, we get that $\kappa_g=\kappa\sin\theta$, $\kappa_n=\kappa\cos\theta$ and $\tau_g\equiv 0$ along $C$. Moreover, since $b=0$, from equation \eqref{EL2}, it directly follows that $H=-c_o$ on $C$.

Next, we use this information in \eqref{EL4}, obtaining 
$$J'\cdot n=\sin\theta\left(2\alpha\kappa''+\left[\alpha\kappa^2-\beta\right]\kappa\right)=0\,.$$
Since $\theta\neq 0$, $\pm\pi$, then $\sin\theta\neq 0$ and $2\alpha\kappa''+(\alpha\kappa^2-\beta)=0$ must hold, which is the equation describing planar elasticae. Therefore, since $C$ is simple and closed, it is a circle of radius $\sqrt{\alpha/\beta}$. In other words, $\alpha\kappa^2=\beta$ holds along $C$. This simplifies equation \eqref{EL3},
$$J'\cdot\nu=0=a\partial_n H\,.$$
Finally, as in previous proposition, the Cauchy-Kovalevskaya Theorem implies that $H\equiv -c_o$ holds on $\Sigma$. {\bf q.e.d.}
\\

Note that if the contact angle is an arbitrary constant, the result may fail. For instance, consider the case of the Willmore energy with elastic boundary $W[\Sigma]$, \eqref{Willmoreelastic} (i.e. $W=E_{a,c_o=0,b=-a,\alpha,\beta}$) and assume that $\theta\equiv 0$ holds on a planar boundary component $C$. Then, the surface $X(\Sigma)$ meets orthogonally the plane containing $C$. Furthermore, since $H=\kappa_n=\kappa$ and $\kappa_g=\tau_g=0$ hold on $C$, it is easy to check that
$$K=\kappa_n\left(2H-\kappa_n\right)-\tau_g^2=\kappa_n^2=H^2$$
also holds on $C$, i.e. $C$ is composed entirely of umbilic points. In this setting, and under the assumption of $\Sigma$ being a topological disc, we conclude from Theorem 7.2 of \cite{P} that $X(\Sigma)$ must be a domain in a sphere of radius $\sqrt{\alpha/\beta}$ whose boundary is the equator. (This domain is the limiting case for the condition $H^2\leq \beta/\alpha$, \cite{Heinz}.)

However, if no restriction on the genus of $\Sigma$ is made, then there exist examples of non totally umbilical Willmore surfaces with non constant mean curvature meeting the plane in a right angle, \cite{BB}.

\section{Global Results of Topological Discs}

Throughout this section, let $X:\Sigma \rightarrow{\bf R}^3$ be an immersion of a topological disc. First, if $\Sigma$ has constant mean curvature, adapting an argument due to Nitsche \cite{N}, we can prove the following result.

\begin{theorem}\label{disc} Let $X:\Sigma \rightarrow {\bf R}^3$ be an immersion of a constant mean curvature surface of disc type having vanishing normal curvature along the boundary. Then, the surface is a planar domain.
\end{theorem}
{\it Proof.\:} We may assume that the surface is given by a conformal immersion of the unit disc $X:D\rightarrow {\bf R}^3$. Let $z$ be the usual complex coordinate in the disc and let $\omega:=\log z$. Although $\omega$ is not well defined, $d\omega=dz/z$ is well defined in $D\setminus \{0\}$. In a neighborhood of $\partial D$, we express the fundamental forms of the immersion as 
\begin{eqnarray*}
&ds_X^2:=e^\rho \lvert d\omega\rvert^2\:,\\
&h:=\frac{1}{2}\Re\left( {\widetilde \Phi }\:d\omega^2+H e^\rho\lvert d\omega\rvert^2\right),
\end{eqnarray*}
where ${\widetilde \Phi}=-\left(H-L_{11}+iL_{12}\right)$, and $L_{ij}$, $i,j=1,2$ are the coefficients of the second fundamental form. The quantities appearing above are known, \cite{Ho}, to satisfy the Gauss and Codazzi equations. In particular, in the case where $H$ is constant, the Codazzi equations reduce to
$${\widetilde \Phi }_{\bar \omega}=0\:.$$
This means that ${\widetilde \Phi}$ defines a holomorphic function. Note that on $\partial D$, we have
\begin{equation}
{\widetilde \Phi } =e^\rho\left(H-i\tau_g\right),\nonumber
\end{equation}
since $\kappa_n\equiv 0$ by hypothesis.

We first assume that $H= 0$ holds. We use the transformation law for quadratic differentials, to obtain the following relation between the Hopf differential in the $\omega$ and $z$ coordinates,
$${\widetilde \Phi }\:d\omega^2={\widetilde \Phi }\omega_z^2\:dz^2={\widetilde \Phi }\left(\frac{1}{z^2}\right)dz^2=:\Phi dz^2\:.$$
In contrast to ${\widetilde \Phi}$, $\Phi$ is globally defined and holomorphic on $D$, as is $z^2\Phi$. The calculation above shows that ${\widetilde \Phi}=z^2\Phi$ on $\partial D$, so we conclude that $\Re\left(z^2\Phi\right)\equiv 0$ holds on $\partial D$. It follows that
on $D$, $\Phi=i c/z^2$ holds for a real constant $c$, which is impossible unless $c=0$, and hence $\Phi$ vanishes identically.

Next, we consider the case where $H$ is a constant different from zero. Consider the image of $\partial D$ under the map $-{\widetilde \Phi}$. Since the real part of $-{\widetilde \Phi}$ never vanishes, this image is contained in a half plane, so  it is clear that the total variation of $\arg {\widetilde \Phi}$ over $\partial D$ vanishes. We write this as
$${\rm Var} \arg {\widetilde \Phi}\:|_{\partial D}=0\:.$$
However, since ${\widetilde \Phi }=z^2\Phi$ on $\partial D$, we get
\begin{eqnarray*}
0&=&{\rm Var} \arg {\widetilde \Phi}\:|_{\partial D}={\rm Var} \arg z^2\Phi\:|_{\partial D}
\\&=&{\rm Var}\arg z^2\:|_{\partial D}+{\rm Var}\arg \Phi \:|_{\partial D}=4\pi+{\rm Var}\arg \Phi \:|_{\partial D}\:.
\end{eqnarray*}
Unless $\Phi\equiv 0$ holds, this gives a contradiction, since, by the Argument Principle, the last term is equal to the total number of zeros of $\Phi$ (which are the umbilics) in $D$, counting multiplicities. So, in particular, the last term is non negative.

Finally, if $\Phi \equiv 0$ holds in $D$, then every point is planar and, hence, the surface is also planar. {\bf q.e.d.}
\\

The previous result simplifies the work for the following classification of constant mean curvature critical discs.

\begin{theorem}\label{clasdisc} Let $X:\Sigma\rightarrow {\bf R}^3$ be a constant mean curvature surface of disc type that is critical for $E_{a,c_o,b,\alpha, \beta}$. Then one of the following occurs:
\begin{enumerate}[(i)]
\item Case $b\neq 0$ and $c_o=0$. Then, the surface is a planar disc or $a=-b>0$ and the surface is a spherical cap satisfying  $H\neq -c_o$. In either case, the boundary of the surface is a circle of radius $\sqrt{\alpha/\beta}$.
\item Case $b=0$. Then $H=-c_o$ and the boundary of the surface is either a circle of radius $\sqrt{\alpha/\beta}$ or a simple closed elastic curve representing a torus knot of type $G(q,1)$ for $q>2$.
\end{enumerate}
In particular, there is no critical constant mean curvature disc when $b\neq 0$ and $c_o \neq 0$.
\end{theorem}
{\it Proof.\:} For constant mean curvature disc type surfaces critical for $E$, we have the following cases, depending on the saddle-splay modulus $b$:
\begin{enumerate}[(i)] 
\item If $b\neq 0$, then by Theorem \ref{disc}, either $H\neq -c_o$ or the surface is planar. In the former case, by Theorem \ref{clas}, $a=-b$ and $c_o=0$ must hold and the surface is a spherical cap bounded by a circle of radius $\sqrt{\alpha/\beta}$. In the latter case, $H=-c_o=0$ and the surface is a disc of radius $\sqrt{\alpha/\beta}$ since the boundary is a simple closed curve critical for the functional $C\mapsto \int_C \left(\alpha\kappa^2+\beta\right)$.
\item If $b=0$, then by \eqref{EL2}, $H\equiv -c_o$ holds. By Proposition \ref{criticality}, the boundary is a simple closed critical curve for $F$, either having constant curvature or representing a $(q,p)$-torus knot, $G(q,p)$. When $b=0$ (i.e. $\mu=0$), critical curves with constant curvature are circles of radii $\sqrt{\alpha/\beta}$. (\emph{We remark that, surprisingly, it is unknown if a circle bounds an embedded constant mean curvature topological disc other than a flat disc or spherical cap.})

Now consider  the case where the boundary curve has non constant curvature and apply Proposition \ref{genus} to get that
$$0=2g\geq \left(p-1\right)\left(q-1\right)\geq \left(p-1\right)\left(2p-1\right),$$
since $q>2p$ (recall that the case $\mu=0$ is the classical elastic energy studied in \cite{ls}) and $p\geq 1$. From this we conclude that $p$ equals one and, hence, the knot is trivial (unknotted). 

In the non constant curvature case, the boundary curve is a critical curve for $F$ representing a torus knot of type $G(q,1)$ for $q>2$. (Examples of minimal discs critical for $E_{a,c_o=0, b=0, \alpha, \beta}$ bounded by elastic curves of type $G(q,1)$ with $q>2$ have been shown in Figure \ref{minelas}.)
\end{enumerate}
We now prove the last statement. If $c_o\neq 0$, then by Theorem \ref{clas}, the mean curvature of a critical constant mean curvature surface must equal $-c_o$. If $b\neq 0$, then $\kappa_n\equiv 0$ on the boundary by \eqref{EL2}. So such a critical constant mean curvature surface cannot exist by Theorem \ref{disc}. {\bf q.e.d.}
\\

Below, we give a lower bound for the energy $E_{a,c_o,b,\alpha, \beta}$ among topological discs that is sharp for some choices of the parameters. For this purpose, we first need to obtain a lower bound for the total squared curvature of a closed curve. The following result is probably well known, we include it for completeness.

\begin{lemma}\label{lemma} If $C(s)$ is any sufficiently smooth closed regular curve in ${\bf R}^3$, then
\begin{equation}\label{lb}
\frac{4\pi^2}{\mathcal{L}[C]}\leq \oint_C \kappa^2\:ds\:,
\end{equation}
with equality holding if and only if the curve is a circle.
\end{lemma}
{\it Proof.\:} We can assume $C(s)=\left(c_1(s), c_2(s), c_3(s)\right)$ is parameterized by arc length. Clearly the integral of $c_i'(s)$ over the curve is zero, so by Wirtinger's inequality,
$$\frac{4\pi}{\mathcal{L}^2}\int_C \left(c_i'(s)\right)^2ds\leq \int_C \left(c_i''(s)\right)^2ds\:,$$
with equality if and only if $c'_i=A_i \cos (2\pi s/\mathcal{L})+B_i\sin(2\pi s/\mathcal{L})$. Summing this over $i$ gives \eqref{lb}. (For simplicity we denote during the proof $\mathcal{L}[C]$ just by $\mathcal{L}$.) 

Next, if ${\vec A}=(A_1,A_2,A_3)$ and ${\vec B}=(B_1,B_2,B_3)$, we get that $C'$ is contained in the span of ${\vec A}$ and ${\vec B}$ so the curve is planar. Also
$$1=\lVert C'(s)\rVert^2=\frac{4\pi}{\mathcal{L}^2}\left(\lvert {\vec A}\rvert^2+\lvert {\vec B}\rvert^2+2\cos(2\pi s/\mathcal{L})\sin(2\pi s/\mathcal{L}){\vec A}\cdot {\vec B}\right),$$
so ${\vec A}$ and ${\vec B}$ are orthogonal. Integrating $C'(s)$ we get 
$$C(s)=\frac{\mathcal{L}}{2\pi}\left(\sin(2\pi s/\mathcal{L}){\vec A}-\cos(2\pi s/\mathcal{L}){\vec B}\right)+{\vec D}\,,$$ 
for a constant vector ${\vec D}$, which is a circle. {\bf q.e.d.}
\\

Using this result, we obtain the following lower bound for some choices of the energy $E$ among topological discs.

\begin{lemma}\label{mindisc} For constants $\alpha>0$ and $\beta>0$, let ${\underline E}:=2\sqrt{\alpha\,\beta}-\lvert b\rvert$ and define the quantity
\begin{equation}\label{ED}
E_D:=2\pi\left({\underline E}+b\right)=2\pi\left(2\sqrt{\alpha\,\beta}-\lvert b\rvert+b\right).\nonumber
\end{equation}
If ${\underline E}\geq 0$ holds, then $E_D$ is a lower bound for the energy $E$ among all sufficiently smooth immersed topological discs in ${\bf R}^3$. For equality to hold, the surface must be a constant mean curvature $H=-c_o\leq 0$ immersion bounded by a circle of radius $\sqrt{\alpha/\beta}$. Moreover, either $b=0$ or $\kappa_n\equiv 0$ must hold along the boundary circle.
\end{lemma}
{\it Proof.\:} Let $\Sigma$ be a topological disc. For any immersion $X:\Sigma\rightarrow{\bf R}^3$, since $a\left(H+c_o\right)^2\geq 0$, the energy $E$ satisfies
\begin{eqnarray*}
E[\Sigma]&=&\int_\Sigma \left(a\left[H+c_o\right]^2+bK\right)d\Sigma+\oint_{\partial \Sigma}\left(\alpha\kappa^2+\beta\right)ds\\
&\ge&b\int_\Sigma K\:d\Sigma+\oint_{\partial \Sigma}\left(\alpha\kappa^2+\beta\right)ds=\oint_{\partial \Sigma}\left(\alpha\kappa^2+b\kappa_g+\beta\right)ds+2\pi b\\
&=&\oint_{\partial \Sigma}\left(\alpha\kappa_n^2+\alpha\kappa_g^2+b\kappa_g+\beta\right)ds+2\pi b\:,
\end{eqnarray*}
where we have used the Gauss-Bonnet Theorem, \eqref{GB}, in the second line. Equality can only hold if $H+c_o\equiv 0$. 

Since the surface is a topological disc, its boundary  is a single closed curve $C$. By the inequality of the arithmetic and geometric means, we get that for all $\epsilon>0$
$$\oint _{C} \lvert b\rvert \, \lvert k_g\rvert\:ds \le \frac{\epsilon \lvert b\rvert}{2}\oint_{C}\kappa_g^2\:ds+\frac{\lvert b\rvert}{2\epsilon} \oint_{C}\:ds\:,$$
with equality if and only if $\kappa_g\equiv{\rm constant}\neq 0$ along $C$ and $\epsilon=1/\lvert \kappa_g\rvert$.

We then obtain for all $\epsilon>0$,
\begin{eqnarray*}
E[\Sigma]-2\pi b&\ge&\oint_{\partial \Sigma}\left(\alpha\kappa_n^2+\alpha\kappa_g^2+b\kappa_g+\beta\right)ds\\
&\geq& \oint_{\partial \Sigma}\left(\alpha\kappa_n^2+\alpha\kappa_g^2-\lvert b\kappa_g\rvert+\beta\right)ds\\
&\ge& \oint_{C}\left(\alpha\kappa_n^2+\alpha\kappa_g^2+\beta\right)ds-\frac{\epsilon \lvert b\rvert}{2}\oint_{C}\kappa_g^2\:ds-\frac{\lvert b\rvert}{2\epsilon} \oint_{C}\:ds\\
&=&\oint_{C}\left(\alpha\kappa_n^2+\left[\alpha- \frac{\epsilon\lvert b\rvert}{2}\right]\kappa_g^2+\left[\beta-\frac{\lvert b\rvert}{2\epsilon}\right]\right)ds\\
&\ge&\oint_{C}\left(\left[\alpha- \frac{\epsilon\lvert b\rvert}{2}\right]\kappa^2+\left[\beta-\frac{\lvert b\rvert}{2\epsilon}\right]\right)ds
\end{eqnarray*}
with equality in the second line if and only if $b\kappa_g=-\lvert b\kappa_g\rvert$. The equality in the last two lines holds if and only if $\kappa_n\equiv 0$ or $b=0$ on $C$. We take $\epsilon=\sqrt{\alpha/\beta}$ to obtain 
\begin{eqnarray*}
E[\Sigma]&\geq& \left(\sqrt{\alpha\,\beta}-\frac{\lvert b\rvert}{2}\right)\oint_{C}\left[\sqrt{\frac{\alpha}{\beta}}\kappa^2+\sqrt{\frac{\beta}{\alpha}}\right]ds+2\pi b\\&=&\frac{1}{2}\,{\underline E}\,\oint_{C}\left[\sqrt{\frac{\alpha}{\beta}}\kappa^2+\sqrt{\frac{\beta}{\alpha}}\right]ds+2\pi b\,.
\end{eqnarray*}
By hypothesis, ${\underline E}\geq 0$, so we can apply Lemma \ref{lemma} to get
$$E[\Sigma]\geq \frac{1}{2}\,{\underline E}\left(\frac{4\pi^2}{\mathcal{L}[C]}\sqrt{\frac{\alpha}{\beta}}+\sqrt{\frac{\beta}{\alpha}}\,\mathcal{L}[C]\right)+2\pi b\,.$$
Recall that equality holds if and only if the boundary component $C$ is a circle. Next, if we minimize the right hand side over all values of $\mathcal{L}[C]$, we find that the minimum occurs at $\mathcal{L}[C]=2\pi\sqrt{\alpha/\beta}$ which yields
$$E[\Sigma]\geq E_D\:.$$
This finishes the proof. {\bf q.e.d.}
\\

We consider first the case $b=0$ and $c_o\neq 0$. It turns out that the existence of constant mean curvature minimizers for $E_{a,c_o,b=0,\alpha,\beta}$ is guaranteed for some choices of the energy parameters.

\begin{prop} If $\alpha$, $\beta$ and $c_o$ are positive constants satisfying $c_o^2\leq \beta/\alpha$ and $b=0$, then there exists a spherical cap with constant mean curvature $H=-c_o$ and this spherical cap is the absolute minimizer of the functional $E_{a,c_o,b=0,\alpha,\beta}$ among all disc type surfaces. 

If $c_o^2>\beta/\alpha$ then no such spherical cap exists.
\end{prop}
{\it Proof.\:} Let $C$ be a circle of radius $\sqrt{\alpha/\beta}$. If $c_o^2\leq \beta/\alpha$ holds, then the sphere of radius $\lvert c_o\rvert^{-1}$ contains the circle $C$ as the boundary of, at least, one spherical cap $\Sigma$. By Lemma \ref{mindisc}, $\Sigma$ minimizes the energy $E_{a,c_o,b=0,\alpha,\beta}$ since,
$$E[\Sigma]=\int_\Sigma a\left(H+c_o\right)^2\,d\Sigma+\oint_{\partial\Sigma}\left(\alpha\kappa^2+\beta\right)ds=4\pi\sqrt{\alpha\,\beta}=E_D\,.$$

The second statement is clear. If $c_o^2>\beta/\alpha$ holds, then $\sqrt{\alpha/\beta}>\lvert c_o\rvert^{-1}$, so a sphere with $H=-c_o$ cannot contain the circle of radius $\sqrt{\alpha/\beta}$. {\bf q.e.d.}
\\

From Theorem \ref{clasdisc}, if $b\neq 0$ and $c_o\neq 0$, it follows that there is no constant mean curvature disc type surface which is a critical point of $E$ and, as a consequence, if ${\underline E}\geq 0$, the lower bound $E_D$ cannot be attained. For these cases, the infimum of the energy is finite and we make the following conjecture.

\begin{conjecture} For positive constants $\alpha$, $\beta$, $c_o$ and $b\neq 0$, if ${\underline E}:=2\sqrt{\alpha\,\beta}-\lvert b\rvert\geq 0$ holds then, the infimum of the energy among all topological discs is attained by an axially symmetric surface with non constant mean curvature.
\end{conjecture}

\begin{remark} In Lemma \ref{minan} below, we show that the infimum of the energy is $-\infty$ in the following cases:
\begin{itemize}
\item $c_o> 0$ and ${\underline E}<0$,
\item $c_o=0$, ${\underline E}<-a$ and $a+b<0$,
\item  $c_o=0$, ${\underline E}<0$ and $b>0$.
\end{itemize}
\end{remark}

The following result discusses all the cases where the spontaneous curvature is zero and the energy has a finite infimum.

\begin{theorem}\label{43} Let ${\underline E}:=2\sqrt{\alpha\,\beta}-\lvert b\rvert$ and $c_o=0$. The finite infima of the energies $E_{a,c_o=0,b,\alpha,\beta}$ among all sufficiently smooth topological discs are as follows:
\begin{enumerate}[(i)]
\item If ${\underline E}\geq -a$ and $a+b<0$, the infimum is 
$$2\pi\left(2\sqrt{\alpha\,\beta}+\left[a+b\right]\right)$$
and it is approached by a limit of spherical caps.
\item If $a+b=0$, the minimum is $4\pi\sqrt{\alpha\,\beta}$ and is attained by a spherical cap bounded by a circle of radius $\sqrt{\alpha/\beta}$.
\item If ${\underline E}\geq 0$ and $b> 0$, the minimum $4\pi\sqrt{\alpha\,\beta}$ is attained by a planar disc bounded by a circle of radius $\sqrt{\alpha/\beta}$.
\item If $b\leq 0$ and $a+b>0$, the minimum $4\pi\sqrt{\alpha\,\beta}$ is attained by a planar disc bounded by a circle of radius $\sqrt{\alpha/\beta}$.
\end{enumerate}
\end{theorem}
{\it Proof.\:} We begin considering the case $a+b\leq 0$ and ${\underline E}\geq -a$ (this is equivalent to ${\underline E}\geq \lvert a+b\rvert-\lvert b\rvert$). Then, we write the energy $E$ with $c_o=0$ as
\begin{eqnarray*}
E[\Sigma]&=&a\int_\Sigma\left(H^2-K\right)d\Sigma+\left(a+b\right)\int_{\Sigma}K\,d\Sigma+\oint_{\partial\Sigma}\left(\alpha\kappa^2+\beta\right)ds\\
&\geq &\left(a+b\right)\int_\Sigma K\,d\Sigma+\oint_{\partial\Sigma}\left(\alpha\kappa^2+\beta\right)ds\\
&=& \oint_{\partial\Sigma}\left(\alpha\kappa^2+\left[a+b\right]\kappa_g+\beta\right)ds+2\pi\left(a+b\right).
\end{eqnarray*}
In the last line, we have used the Gauss-Bonnet Theorem \eqref{GB}. For the inequality, we have used that $H^2-K\geq 0$ where equality holds if and only if the surface is totally umbilical. Next, we follow the exact same steps as in the proof of Lemma \ref{mindisc} with $b$ replaced by $a+b$ to conclude
$$E[\Sigma]\geq 2\pi\left(2\sqrt{\alpha\,\beta}+\left[a+b\right]\right).$$ 
Here we can apply Lemma \ref{lemma} since ${\underline E}\geq \lvert a+b\rvert-\lvert b\rvert$ holds. Recall that equality above holds if and only if the boundary is a circle of radius $\sqrt{\alpha/\beta}$ (and the surface is totally umbilical).  

Assume now that $a+b<0$ and consider $\Sigma_R$ to be a sequence of spherical caps of radii $R$ bounded by a circle of radius $\sqrt{\alpha/\beta}$. We then have
$$\oint_{\partial\Sigma_R}\kappa_n^2\,ds=\frac{2\pi}{R}$$
and so, $E[\Sigma_R]\longrightarrow 2\pi\left(2\sqrt{\alpha\,\beta}+\left[a+b\right]\right)$. No minimizer can exist since it would necessarily be a totally umbilical surface bounded by a circle of radius $\sqrt{\alpha/\beta}$ on which $\kappa_n\equiv 0$ holds and no such circle exists.

For part (ii), let $a+b=0$. Then, the condition ${\underline E}\geq -a$ is automatically satisfied and previous inequalities hold. Let $\Sigma$ be a spherical cap bounded by a circle of radius $\sqrt{\alpha/\beta}$. Note that the radius of the sphere must be bigger than or equal $\sqrt{\alpha/\beta}$. Then, all equalities in above estimate hold and $\Sigma$ attains the minimum $4\pi\sqrt{\alpha/\beta}$. (See Theorem \ref{clas}.)

Case (iii) follows immediately from Lemma \ref{mindisc}. Equalities in that lemma for our parameters hold if and only if the surface is minimal bounded by a circle of radius $\sqrt{\alpha/\beta}$ satisfying $\kappa_n\equiv 0$ along it. Therefore, it is the planar critical disc and it attains the minimum $E_D=4\pi\sqrt{\alpha\,\beta}$.

For the last case, assume $a+b>0$ and $b\leq 0$. Then, we write the energy $E$ as
\begin{eqnarray*}
E[\Sigma]&=&\left(a+b\right)\int_\Sigma H^2\,d\Sigma-b\int_{\Sigma}\left(H^2-K\right)d\Sigma+\oint_{\partial\Sigma}\left(\alpha\kappa^2+\beta\right)ds\\
&\geq &\oint_{\partial\Sigma}\left(\alpha\kappa^2+\beta\right)ds\geq 4\pi\sqrt{\alpha\,\beta}\,,
\end{eqnarray*}
arguing as above. Equalities hold if and only if the surface is minimal and totally umbilical, i.e. planar, and it is bounded by a circle of radius $\sqrt{\alpha/\beta}$. Therefore, the critical planar disc attains the minimum $4\pi\sqrt{\alpha\,\beta}$. {\bf q.e.d.}
\\

To summarize, in this section, we have characterized those cases in which the infimum of the energy among topological discs is finite. In each of these cases, we have either produced a constant mean curvature minimizer, produced a constant mean curvature minimizing sequence or we have shown that no constant mean curvature minimizer can exist.

\section{Global Results of Topological Annuli}

In this section we investigate critical annuli. Surprisingly, the results are more complete than in the case of topological discs.

Let $X:\Sigma\rightarrow{\bf R}^3$ be the immersion of a topological annulus $\Sigma$ with constant mean curvature $H=-c_o$ critical for $E$. Then, from Proposition \ref{criticality}, both boundary components are critical for $F$ and represent torus knots, namely, $G_1$ and $G_2$. By a result from knot cobordism theory \cite{Lither}, both torus knots are of the same type. Some examples of critical annuli for $E_{a,c_o=0,b=0,\alpha,\beta}$ have been shown in Figure \ref{minelasan}. In particular, in Figure \ref{minelasan} (C) the boundary components are non trivial knots.

However, the simplest case occurs when both boundary components are circles as discussed in the following result.

\begin{prop}\label{rotational} Let $X:\Sigma\rightarrow {\bf R}^3$ be an immersion of a constant mean curvature compact surface with boundary $\partial\Sigma$. If any connected component of $\partial\Sigma$ is a circle on which $\tau_g\equiv 0$ holds, then the surface is axially symmetric.  
\end{prop}
{\it Proof.\:} Let $C$ be a circle representing a connected component of $\partial\Sigma$. We can assume that the circle $C$ lies in a horizontal plane. Let $\mathcal{R}_t$ denote the one parameter family of rotations about a vertical axis passing through the circle's center. For any surface $\Sigma$ with constant mean curvature $H$, the function 
$$\psi:=\partial_t \left(\mathcal{R}_tX\right)_{t=0}\cdot \nu=E_3\times X\cdot \nu$$ 
is the normal part of the derivative of a variation of $X$ through constant mean curvature surfaces and, as such, $\psi$ defines a Jacobi field on $\Sigma$, that is $L[\psi]:=\Delta \psi +\lVert d\nu\rVert ^2\psi=0$ holds. Since $C$ is invariant under $\mathcal{R}_t$, we have that $\psi\equiv 0$ along $C$. Also, along $C$, we have
\begin{eqnarray*}
\partial_n\psi& =&E_3\times n\cdot \nu +E_3\times X\cdot d\nu(n)\\
&=&-T\cdot E_3-\tau_g E_3\times X\cdot T -\left(2H-\kappa_n\right)E_3\times X\cdot n\\
&=&0\:,
\end{eqnarray*} 
since $E_3$ is normal to $C$, $\tau_g\equiv 0$ and $E_3\times X$ is tangent to $C$.

We now use the well known fact that constant mean curvature surfaces are real analytic. It follows that equation $L=0$ is a second order elliptic PDE with $\mathcal{C}^\omega$ coefficients. The circle $C$ is a real analytic curve which we view as the initial curve for the Cauchy problem $Lf=0$ with $f=\partial_nf=0$ on $C$. By the Cauchy-Kovalevskaya Theorem, this problem, locally, has a unique analytic solution $f\equiv 0$, so we obtain $\psi\equiv 0$ locally, since every solution of $L=0$ is analytic by elliptic regularity. Again, using analyticity of the surface, we obtain that $\psi\equiv 0$ {\it globally} on $\Sigma$ and it follows that $X(\Sigma)$ is axially symmetric. {\bf q.e.d.} 
\\

The axially symmetric constant mean curvature surfaces are known as Delaunay surfaces, \cite{D}. They fall into six types: planes, spheres, right circular cylinders, catenoids, unduloids and nodoids. Identifying ${\bf R}^2$ with the complex plane ${\bf C}$, these surfaces can be represented as $X(u,\vartheta)=(r(u)e^{i\vartheta} ,z(u))$, where $r$ and $z$ satisfy equations
\begin{eqnarray}
&ur+Hr^2\equiv {\rm constant}:=\varpi\,,\label{nodc}\\
&dz=r_u\:dw\:.\label{nodc2}
\end{eqnarray}
Here $\nu=(ue^{i\vartheta}, w)$ is the Gauss map of $X$, $H$ is the constant mean curvature which we assume, after our choice of orientation, to always be non positive and $\varpi$ is a constant which we call the \emph{flux parameter}. In \eqref{nodc2} it is understood that the quadratic in \eqref{nodc} has been solved for $r=r(u)$. 

For  minimal surfaces, i.e. $H=0$, we have: for $\varpi=0$, the surface is a plane, while for $\varpi\neq 0$, the surface is a catenoid. Next, if $H<0$, for $\varpi=0$ we have either a right circular cylinder or a sphere (depending on $u$ being constant or not, respectively), for $\varpi<0$ the surface is a nodoid, while for $\varpi>0$ with $1+4\varpi H\geq 0$ we get an unduloid.

The first result concerns critica when the saddle-splay modulus $b$ is zero.

\begin{prop} In a Delaunay surface with $H=-c_o$, any embedded annular domain $\Omega$ bounded by two parallels of radii $\sqrt{\alpha/\beta}$ is critical for $E$ with $b=0$. Moreover, the energy of any of these domains is
$$E[\Omega]=8\pi\sqrt{\alpha\,\beta}\,.$$
\end{prop}
{\it Proof.\:} Consider a Delaunay surface with $H=-c_o$ and let $\Omega$ be an annular domain bounded by two parallels of radii $\sqrt{\alpha/\beta}$. In the case of a nodoid, we assume that the surface is embedded.

Since $H=-c_o$ and $b=0$, equations \eqref{EL1} and \eqref{EL2} are clearly satisfied. Moreover, using that both boundary components are circles of radii $\sqrt{\alpha/\beta}$ (i.e. $\alpha\kappa^2=\beta$) and that $\tau_g\equiv 0$ holds along them, it is easy to check that equations \eqref{EL3}-\eqref{EL4} are also satisfied, proving that the domains $\Omega$ are critical for $E$ with $b=0$.

Finally, we compute their energy. Once more, since $H=-c_o$ and $b=0$, we obtain
\begin{eqnarray*}
E[\Omega]&=&\int_\Omega \left(a\left[H+c_o\right]^2+bK\right)d\Omega+\oint_{\partial\Omega}\left(\alpha\kappa^2+\beta\right)ds\\
&=&\oint_{\partial\Omega}\left(\alpha\kappa^2+\beta\right)ds=2\beta\mathcal{L}[\partial\Omega]=8\pi\sqrt{\alpha\,\beta}\,.
\end{eqnarray*}
The equalities in the last line follow since the boundary is composed by two circles of radii $\sqrt{\alpha/\beta}$ (i.e. of curvature $\sqrt{\beta/\alpha}$). {\bf q.e.d.}
\\

We consider now the case $b\neq 0$. In this case, when $H=-c_o$, by equation \eqref{EL2}, $\kappa_n\equiv 0$ must hold along the boundary. Since the normal curvature of parallels is given by $\kappa_n=-u/r$, this means that, except possibly in planes and nodoids, the annular domains of Delaunay surfaces bounded by two parallels are never critical. In contrast, there are infinitely many annular domains in a nodoid satisfying $\kappa_n\equiv 0$ on the boundary.

\begin{prop}\label{UD} For $a>0$, $b\neq 0$, $c_o>0$, $\alpha>0$ and $\beta>0$ there is a unique nodoid with $H=-c_o$ having critical domains for the functional $E$. This nodoid contains, up to rigid motions, exactly four such critical domains.
\end{prop}
{\it Proof.\:} It is clear that for a surface with $H=-c_o<0$, \eqref{EL1} holds. Since $b\neq 0$, it then follows that $\kappa_n\equiv 0$ must hold on $\partial \Sigma$ which means that $u\equiv 0$ on $\partial \Sigma$. Solving \eqref{nodc} gives 
\begin{equation}\label{r}
r=\frac{u\pm \sqrt{u^2-4\varpi c_o}}{2c_o}\,. 
\end{equation} 
For nodoids, only the plus sign is used in the formula \eqref{r} so setting $u=0$ gives $r=\sqrt{-\varpi/{c_o}}$ on the boundary. Since the geodesic torsion is zero and $\kappa_n$ and $\kappa_g$ are constant on parallels one sees that equation \eqref{ELb} holds if and only if $\alpha\kappa^2=\beta$ holds on the boundary, i.e. this means that both boundary radii are $\sqrt{\alpha/\beta}$. Setting this equal to $\sqrt{-\varpi/{c_o}}$ determines the value of $\varpi$.

Geometrically, is clear that up to rigid motions, there are exactly four embedded domains bounded by circles of radii $\sqrt{\alpha/\beta}$ in the nodoid and on the circles the normal is vertical, i.e. $\kappa_n\equiv 0$ holds. {\bf q.e.d.}
\\

\begin{figure}[h!]
\makebox[\textwidth][c]{
\begin{subfigure}[b]{0.323\linewidth}
\includegraphics[width=\linewidth]{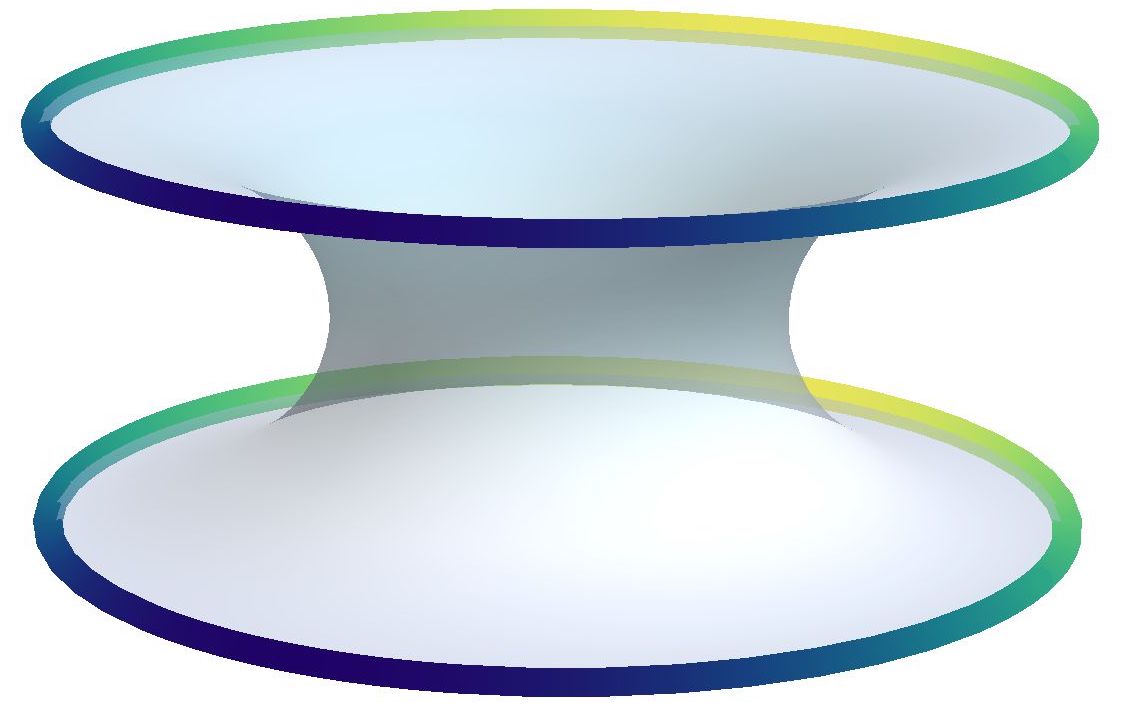}
\caption{$\mathcal{N}_{1}$}
\end{subfigure}
\,
\begin{subfigure}[b]{0.288\linewidth}
\includegraphics[width=\linewidth]{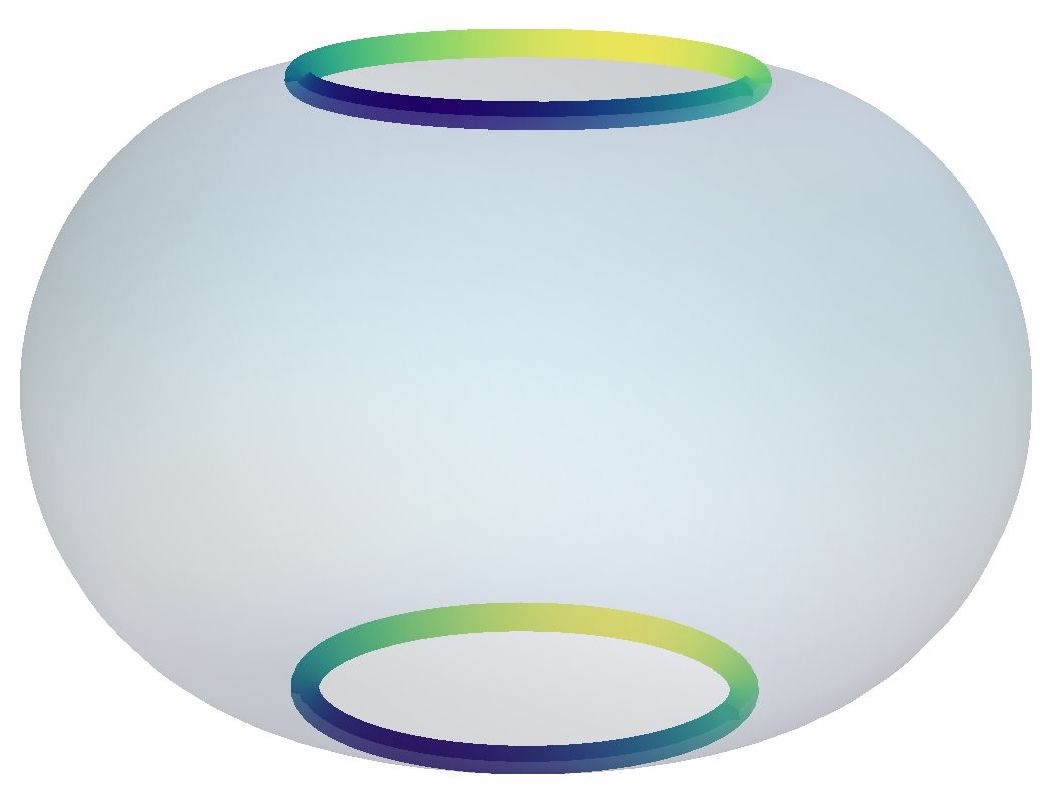}
\caption{$\mathcal{N}_{2}$}
\end{subfigure}
\,
\begin{subfigure}[b]{0.318\linewidth}
\includegraphics[width=\linewidth]{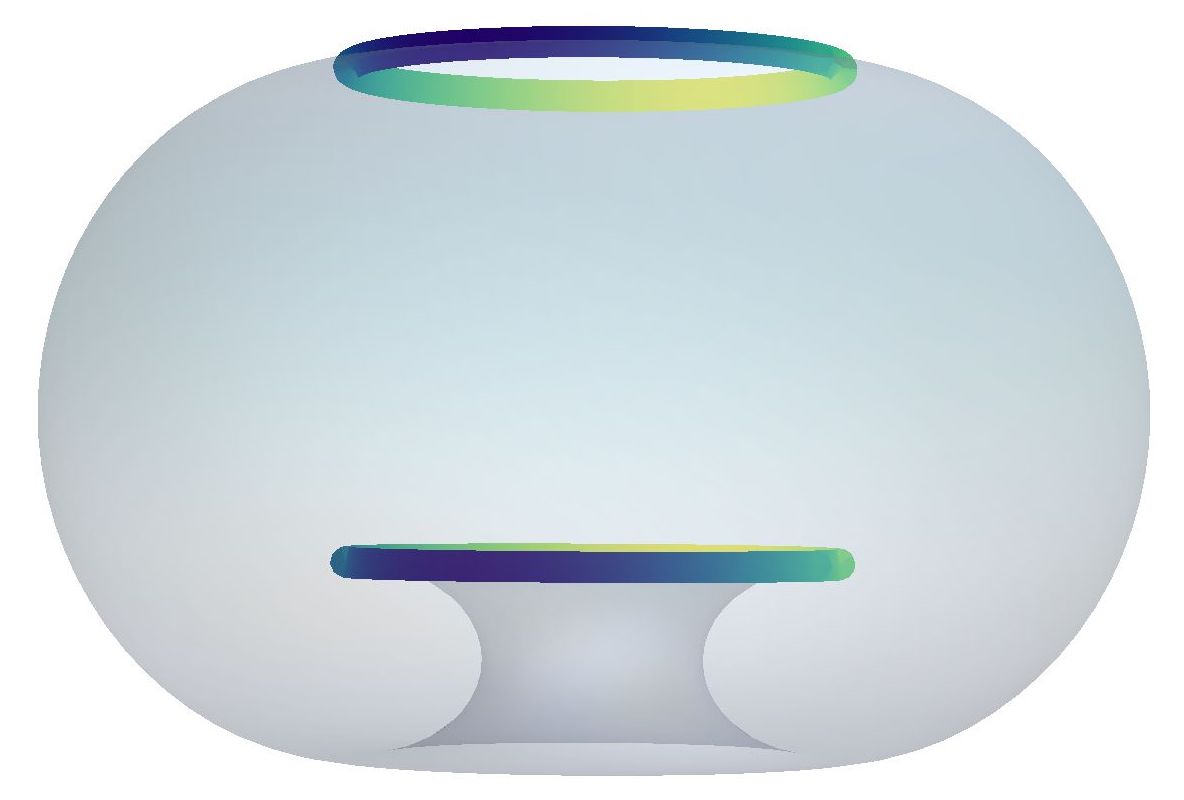}
\caption{$\mathcal{N}_{3}$}
\end{subfigure}
\,
\begin{subfigure}[b]{0.31\linewidth}
\includegraphics[width=\linewidth]{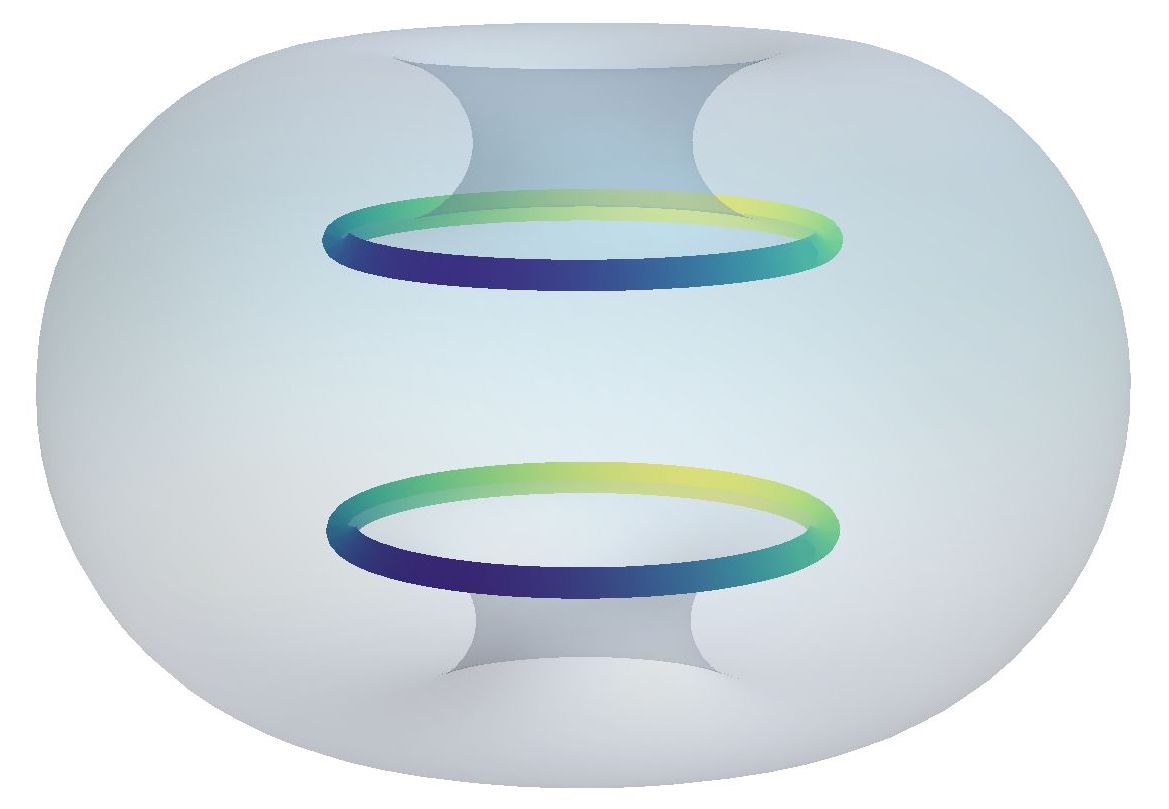}
\caption{$\mathcal{N}_{4}$}
\end{subfigure}}
\caption{The four distinct types of nodoidal domains given in Proposition \ref{UD} critical for $E_{a,c_o>0,b,\alpha,\beta}$.}
\label{nods}
\end{figure}

As a consequence of Proposition \ref{cca} and Proposition \ref{rotational}, above nodoidal domains can be characterized in terms of the contact angle.

\begin{cor} Let $X:\Sigma\rightarrow{\bf R}^3$ be a critical immersion for $E$ with $b\neq 0$ and $c_o\neq 0$. If, at least, one boundary component is planar and the contact angle satisfies $\theta\equiv \pm\pi/2$ along it, then the surface is one of the four nodoidal domains $\mathcal{N}_i$, $i=1,...,4$ with $H=-c_o<0$.
\end{cor}
{\it Proof.\:} Let $C$ be the planar boundary component. Since $\theta\equiv \pm\pi/2$, from \eqref{kg}-\eqref{taug}, $\kappa_g=\pm\kappa$, $\kappa_n\equiv 0$ and $\tau_g\equiv 0$ hold along $C$. Furthermore, applying Proposition \ref{cca}, we have that $X(\Sigma)$ has constant mean curvature $H=-c_o$. In this situation, the Euler-Lagrange equation \eqref{EL4} reduces to
$$J'\cdot n+a\left(H+c_o\right)^2+bK=2\alpha\kappa''+\left(\alpha\kappa^2-\beta\right)\kappa=0\,,$$
up to a sign, which implies that $C$ is a circle of radius $\sqrt{\alpha/\beta}$. 

Finally, since $X(\Sigma)$ has constant mean curvature and $C$ is a circle along which $\tau_g\equiv 0$ holds, it follows from Proposition \ref{rotational} that $X(\Sigma)$ is axially symmetric. Moreover, since $b\neq 0$ and $c_o\neq 0$, from Theorem \ref{clasdisc} we get that $\Sigma$ is an annulus. An analysis of these annular surfaces as done above concludes the proof. {\bf q.e.d.}
\\

To finish this section we are going to analyze the problem of minimizing the energy $E$ among topological annuli. We begin by obtaining lower bounds of the energy in all cases.

\begin{lemma}\label{minan} For constants $\alpha>0$ and $\beta>0$, let ${\underline E}:=2\sqrt{\alpha\,\beta}-\lvert b\rvert$ and define the quantity
\begin{equation}\label{EA}
E_A:=4\pi\,{\underline E}=4\pi\left(2\sqrt{\alpha\,\beta}-\lvert b\rvert\right).\nonumber
\end{equation}
Then, we have:
\begin{enumerate}[(i)]
\item If ${\underline E}\geq 0$ holds, then $E_A$ is a lower bound for the energy $E$ among all sufficiently smooth topological annuli in ${\bf R}^3$.
\item If ${\underline E}<0$ and $c_o>0$ hold, then the infimum of the energy among all topological annuli is $-\infty$.
\item If ${\underline E}<-a$, $c_o=0$ and $a+b<0$ hold, then the infimum of the energy among all topological annuli is $-\infty$.
\item If ${\underline E}< 0$, $c_o=0$ and $b>0$ hold, then the infimum of the energy among all topological annuli is $-\infty$.
\end{enumerate}
The result of cases (ii)-(iv) also holds for topological discs.
\end{lemma}
{\it Proof.\:} We begin proving part (i). Let $X:\Sigma \rightarrow{\bf R}^3$ be the immersion of a topological annulus. Arguing as in the proof of Lemma \ref{mindisc}, for any immersion the energy $E$ satisfies
$$E[\Sigma]\geq \oint_{\partial\Sigma}\left(\alpha\kappa_n^2+\alpha\kappa_g^2+b\kappa_g+\beta\right)ds\,.$$
Here, we have used the Gauss-Bonnet Theorem for annuli, \eqref{GB}.

We write $\partial\Sigma\equiv C_1\cup C_2$ and using the inequalities  on each boundary component $C_i$, $i=1,2$ obtained in Lemma \ref{mindisc}, we get that
$$E[\Sigma]\geq E_A\,.$$

To prove (ii), we first assume that $b<0$ holds. Consider the critical nodoidal domain $\mathcal{N}_2$ with $H=-c_o$. For this surface we have
$$E[\mathcal{N}_2]=\int_\Sigma\left(a\left[H+c_o\right]^2+bK\right)d\Sigma+\oint_{\partial\Sigma}\left(\alpha\kappa^2+\beta\right)ds=E_A\,,$$
since the total curvature of $\mathcal{N}_2$  is $4\pi$ (see Figure \ref{nods}). Also, note that by hypothesis $E_A=4\pi{\underline E}<0$ so $E[\mathcal{N}_2]<0$ holds. The surface ${\mathcal N}_2$ will be represented as the image of an infinite covering
$$X:[\sigma_1,\sigma_2]\times (-\infty, \infty)\rightarrow {\bf R}^3,\quad (\sigma,t)\mapsto (r(\sigma)e^{it}, z(\sigma))\:$$
It is well known that $X$ can be embedded in a $2\pi$ periodic family of isometric constant mean curvature immersions $X_\vartheta:[\sigma_1,\sigma_2]\times (-\infty, \infty)\rightarrow {\bf R}^3$ with $X_0\equiv X$ and $X_\pi$ being an immersion of part of an unduloid. The immersions for $\vartheta\in (0,\pi)$ are invariant under a helicoidal motion. These will all have constant mean curvature $H=-c_o$ and the curves $t\mapsto X_\vartheta(\sigma_i, t)$ will have normal curvature $\kappa_n(\vartheta)\equiv H(1-\cos\vartheta)$. Hence the curvatures of these curves satisfy $[\kappa(\vartheta)]^2=\kappa_g^2+(H[1-\cos\vartheta])^2$. It follows that for a fundamental domain $S_\vartheta$ of $X_\vartheta$, we have
\begin{eqnarray*}
E[S_\vartheta]&=&\int_{S_\vartheta} \left(a\left[H+c_o\right]^2+bK\right)d\Sigma+\sum_{i=1}^{2}\oint_{\sigma=\sigma_i} \left(\alpha \kappa^2+\beta\right)ds\\
&=&E_A+\sum_{i=1}^{2}\left(H\left[1-\cos\vartheta\right]\right)^2\oint_{\sigma=\sigma_i}ds<0
\end{eqnarray*}
for $\vartheta \approx 0$.

Let $\Sigma_m$ denote the surface obtained by taking the m-fold covering of $S_\vartheta$, $X_\vartheta:[\sigma_1,\sigma_2]\times (0, 2\pi m)\rightarrow {\bf R}^3$. If we smooth out the corners of this surface, we obtain that for a positive constant $P$ independent of $m$,
$$E[\Sigma_m]=m\left(E_A+\sum_{i=1}^{2}\left(H\left[1-\cos\vartheta\right]\right)^2\oint_{\sigma=\sigma_i}ds\right)+P \longrightarrow -\infty\:$$
as $m\rightarrow \infty$. Here $P$ represents the integral $\int (\alpha \kappa^2+\beta)\:ds$ over the arcs connecting the helices, the corners at the intersections of these arcs with the helices being suitably smoothed. This shows that the infimum of the energy among topological discs is $-\infty$. To obtain the result for annuli, simply cut a small disc out of $\Sigma_m$. This will add a small positive term, independent of $m$, to the energy, so the limit is not affected.

The case $b>0$ is proved in a similar way with the critical nodoidal domain ${\mathcal N}_1$ replacing ${\mathcal N}_2$.

The proof of case (iii) is also similar to the proof of (ii). Let $\epsilon$ be such that $4\pi>\epsilon>0$ and
$$8\pi\sqrt{\alpha\,\beta}+4\pi\left(a+b\right)< 3 \epsilon \left(a+b\right).$$
For $u>0$ and fixed $H<0$ the formulas (\ref{nodc}) and (\ref{nodc2}) give immersions of a family of positively curved nodoids parameterized by ${\varpi}<0$. These immersions, as functions of $u$ and $\vartheta$, are exactly the inverses of the Gauss maps of the surfaces. We fix a symmetric domain $\Omega \subset S^2$, bounded by circles of the same radii, such that the total curvature of the domain $\Omega$ is larger than $4\pi-\epsilon$ for all ${\varpi}<0$, i.e.
$$\int_\Omega K\:d\Sigma>4\pi-\epsilon\,.$$
This is possible since the curvature integral is just the spherical area. Then, for $\varpi \approx 0$, we have 
$$a\int_\Omega \left(H^2-K\right)d\Sigma <-\epsilon\left(a+b\right)$$ 
so, combining both the last two inequalities, we get
\begin{eqnarray*}
a\int_\Omega H^2\:d\Sigma +b\int_\Omega K\:d\Sigma &=&a\int_\Omega \left(H^2-K\right)d\Sigma +\left(a+b\right)\int_\Omega K\:d\Sigma\\
&<& 4\pi\left(a+b\right)-2\epsilon\left(a+b\right).
\end{eqnarray*}
Note that the integrals given above are scale invariant. If we rescale the nodoids so that both boundary circles have radii $\sqrt{\alpha/\beta}$ then for $\varpi\approx 0$ we have
$$E[\Omega]=8\pi \sqrt{\alpha\,\beta}+4\pi\left(a+b\right)-2\epsilon\left(a+b\right)<\epsilon\left(a+b\right)<0.$$ 

We now proceed as above and consider the isometric deformation of the nodoids through helicoidal surfaces. If $S_\vartheta$ denotes the fundamental domain for one of these helicoidal surfaces, then the energy contribution of its two helical boundary components increments continuously, while the energies of the arcs connecting the helical arcs (sufficiently smoothed at the corners) will add a fixed constant $P$ to the energy. Thus, the energy of $S_\vartheta$, for $\vartheta\approx 0$ satisfies 
$$E[S_\vartheta]=8\pi \sqrt{\alpha\,\beta}+4\pi\left(a+b\right)-2\epsilon\left(a+b\right)+P$$ 
and the energy of $m$ contiguous copies of the fundamental domain, $\Sigma_m$, will satisfy 
$$E[\Sigma_m]<m\left(8\pi \sqrt{\alpha\, \beta}+4\pi\left[a+b\right]-2\epsilon\left[a+b\right]\right)+P\longrightarrow -\infty,$$ 
as $m\rightarrow \infty$. This shows that for discs the energy is not bounded below. The result for annuli follows by cutting out a small disc as before.

Finally, we prove case (iv). Again, this case is similar to the previous cases. Let $S$ be a vertical catenoid which is symmetric with respect to the plane $z=0$. Choose $\epsilon$, $0<\epsilon<4\pi$ such that $8\pi\sqrt{\alpha\,\beta}-b\left(4\pi-\epsilon\right)<0$. There exists $z_o>0$ such that $S_o:=S\cap \{\lvert z\rvert<\lvert z_o\rvert\}$ has total curvature less than $-4\pi+\epsilon$. This inequality is unaffected by rescaling, so we can assume that $S_o$ has boundary circles of radii $\sqrt{\alpha/\beta}$ which gives $E[S_o]<0$. We then proceed as above using domains in the isometric helicoid. {\bf q.e.d.}
\\

In the following theorem we state the infima of the energy $E$ among topological annuli in every case.

\begin{theorem} Assume $\alpha>0$, $\beta>0$, $c_o\geq 0$ and denote ${\underline E}:=2\sqrt{\alpha\,\beta}-\lvert b\rvert$. Then, the infima of the energies $E_{a,c_o\geq 0,b,\alpha,\beta}$ among all sufficiently smooth annuli are given in the following table:
\begin{table}[H]
\makebox[\textwidth][c]{
\begin{tabular}{|c|c|c|c|c|} 
\hline
Case & Parameters & Energy & Surface & Type \\
\hline
(i) & ${\underline E}\ge 0$, $c_o>0$, $b>0$ & $4\pi\left(2\sqrt{\alpha \beta}- b\right)$ & $\mathcal{N}_{1}$ and $\mathcal{N}_{4}$ & Minimum\\ 
(ii) & $c_o>0$, $b=0$ & $8\pi\sqrt{\alpha \beta}$ & Multiple solutions & Minimum\\ 
(iii) & ${\underline E}\ge 0$, $c_o>0$, $b<0$ & $4\pi\left(2\sqrt{\alpha \beta}+b\right)$ & $\mathcal{N}_{2}$ & Minimum\\
(iv) & ${\underline E}\ge 0$, $c_o=0$, $b>0$ & $4\pi\left(2\sqrt{\alpha \beta}- b\right)$ & Limit of catenoid domains & Infimum\\ 
(v) & $c_o=0$, $b=0$ & $8\pi\sqrt{\alpha \beta}$ & Multiple solutions & Minimum\\ 
(vi) & ${\underline E}\ge -a$, $c_o=0$, $a+b<0$ & $4\pi\left(2\sqrt{\alpha \beta}+[a+ b]\right)$ & Limit of spherical annuli & Infimum\\ 
(vii) & $c_o=0$, $a+b=0$ & $8\pi\sqrt{\alpha \beta}$ & Spherical annulus & Minimum\\ 
(viii) & $c_o=0$, $b<0$, $a+b>0$ & $8\pi\sqrt{\alpha \beta}$ & Limit of planar annuli & Infimum\\    
\hline
\end{tabular}}
\end{table}
\noindent In all cases not listed in the table above, the energy is not bounded below.
\end{theorem}
{\it Proof.\:} The surfaces $\mathcal{N}_{i}$, $i=1,...,4$ satisfy $H=-c_o$ and have boundary circles of radii $\sqrt{\alpha/ \beta}$. Moreover, observe in Figure \ref{nods} that the total curvature of $\mathcal{N}_{1}$ and $\mathcal{N}_{4}$ is $-4\pi$, while the total curvature of $\mathcal{N}_{2}$ is $4\pi$. It is thus clear that these surfaces achieve the lower bound of Lemma \ref{minan} for $b>0$ and $b<0$, respectively, and cases (i) and (iii) follow. 

\begin{figure}[h!]
\centering
\begin{subfigure}[b]{0.25\linewidth}
\includegraphics[width=\linewidth]{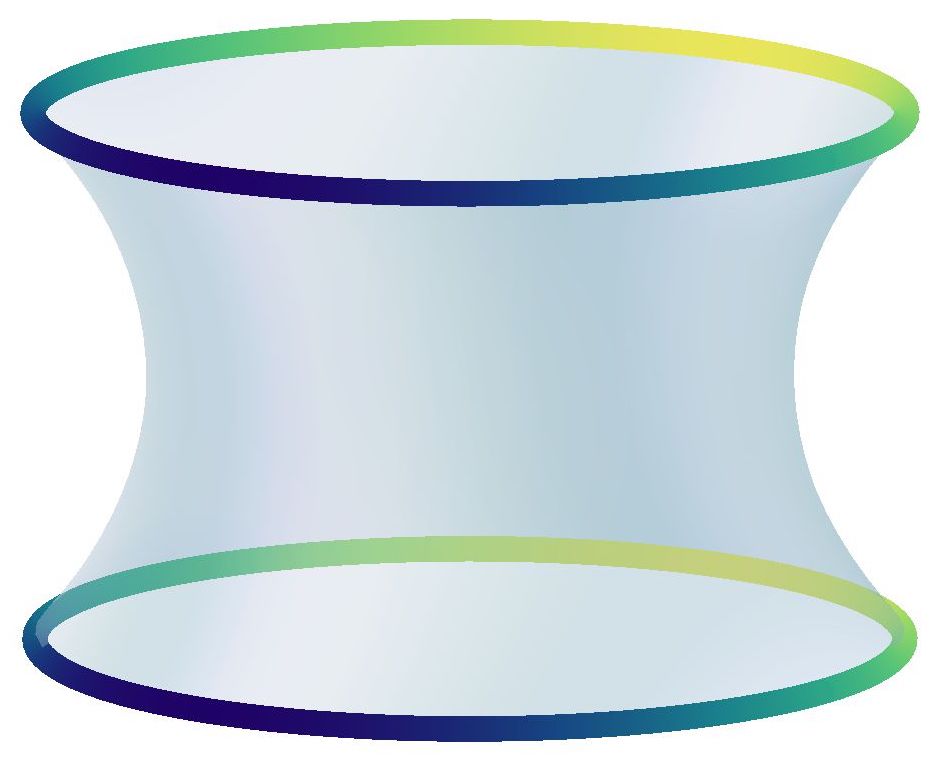}
\end{subfigure}
\,\,\,
\begin{subfigure}[b]{0.31\linewidth}
\includegraphics[width=\linewidth]{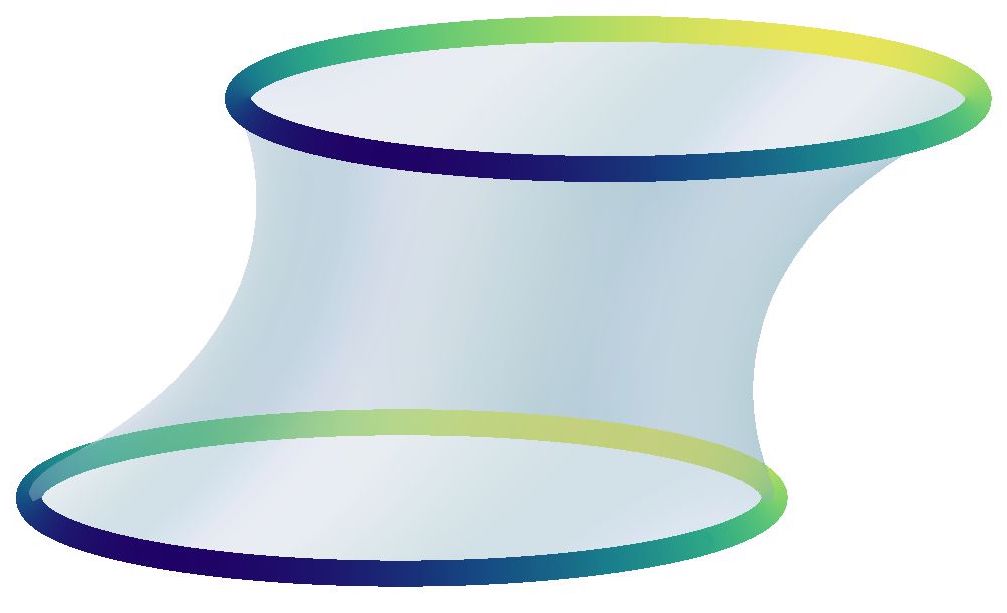}
\end{subfigure}
\begin{subfigure}[b]{0.335\linewidth}
\includegraphics[width=\linewidth]{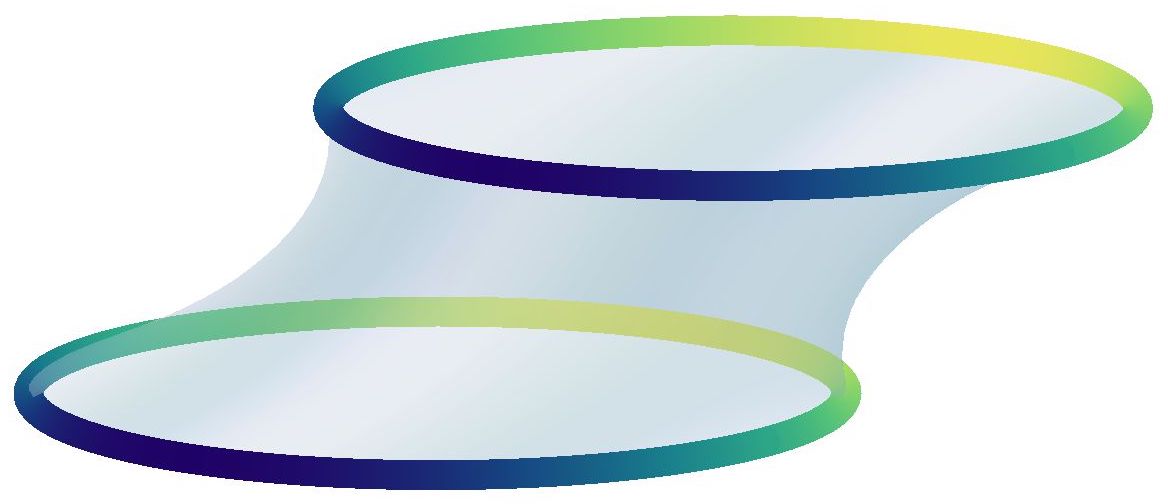}
\end{subfigure}
\caption{A part of the catenoid (left) and two Riemann's minimal examples (center and left) bounded by circles of radii $\sqrt{\alpha/\beta}$. These domains are critical absolute minimizers for $E$ with $b=0$ and $c_o=0$, i.e. for $E_{a,c_o=0,b=0,\alpha,\beta}$.}\label{minimalexamples}
\end{figure}

In case (ii), any domain in a Delaunay surface with $H=-c_o<0$ bounded by two parallels having radii $\sqrt{\alpha/\beta}$ will achieve the lower bound in Lemma \ref{minan}. Note that in this case the restriction ${\underline E}\geq 0$ trivially holds. Similarly, when $c_o=0$ and $b=0$, i.e. case (v), any minimal surface bounded by two circles of radii $\sqrt{\alpha/\beta}$ attains the lower bound in Lemma \ref{minan}. An example of this is a suitable annular domain in a catenoid bounded by two parallels of radii $\sqrt{\alpha/\beta}$. (See Figure \ref{minimalexamples}.)

For case (iv), fix a vertical catenoid $S$ which is symmetric with respect to the plane $z=0$. Let $S_R$ be the region of $S$ bounded by the planes $z=\pm R$ and let $\Sigma_R$ be the domains obtained by rescaling $S_R$ so that the boundary circles have radii $\sqrt{\alpha/\beta}$. It is then clear that
$$\int_{S_R} K\:d\Sigma \longrightarrow -4\pi$$ 
as $R\rightarrow \infty$ and so $E_{a,c_o=0,b>0,\alpha,\beta}[\Sigma_R]\longrightarrow 4\pi\left(2\sqrt{\alpha \beta}-b\right)$. In this case, there is no minimizer. If one were to exist then by Proposition \ref{rotational} and the proof of Lemma \ref{minan}, it would necessarily be a catenoid domain bounded by circles on which $\kappa_n\equiv 0$ holds and there is no such domain.

For cases (vi) and (vii), we argue as in the proof of Theorem \ref{43} to obtain ($E=E_{a,c_o=0,b<0,\alpha,\beta}$)
$$E[\Sigma]\geq 4\pi\left(2\sqrt{\alpha\,\beta}+\left[a+b\right]\right)$$
where equality holds if and only if the surface is totally umbilical bounded by two circles of radii $\sqrt{\alpha/\beta}$ and either $a+b=0$ or $\kappa_n\equiv 0$ holds along the boundary. In case (vii), since $a+b=0$, a domain in a sphere of radius $R>\sqrt{\alpha/\beta}$ bounded by two circles of radii $\sqrt{\alpha/\beta}$ will realize the minimum value (see Theorem \ref{clas}). 

For (vi), i.e. $a+b<0$, we consider $\Sigma_R$ to be the sequence of spherical annuli of radii $R$ bounded by two circles of radii $\sqrt{\alpha/\beta}$. By the same argument as in the proof of Theorem \ref{43}, we get $E[\Sigma_R]\longrightarrow 4\pi\left(2\sqrt{\alpha\,\beta}+\left[a+b\right]\right)$ as $R\rightarrow{\infty}$, but again no minimizer can exist.

Finally, the inequalities needed for case (viii) also follow by the same reasoning as in Theorem \ref{43}. For the annular case, we obtain ($E=E_{a,c_o=0,b<0,\alpha,\beta}$)
$$E[\Sigma]\geq 8\pi\sqrt{\alpha\,\beta}\,,$$
where equality holds if and only if the surface is planar and bounded by two circles of radii $\sqrt{\alpha/\beta}$. Clearly, such a domain in the plane does not exist. However, we can consider the sequence of planar annuli $\Sigma_R$ bounded by two circles of radii $\sqrt{\alpha/\beta}\pm1/R$, respectively. Then, $E[\Sigma_R]\longrightarrow 8\pi\sqrt{\alpha\,\beta}$ when $R\rightarrow\infty$.

The last statement of the theorem follows from Lemma \ref{minan}.  {\bf q.e.d.}
\\

In the case $b=0$ there are multiple domains attaining the minimum of the energy $E_{a,c_o\geq 0,b=0,\alpha,\beta}$. Indeed, we cannot even conclude that an energy minimizing surface is axially symmetric. For $c_o=0$, any annular domain in a Riemann's minimal example bounded by circles of radii $\sqrt{\alpha/\beta}$ also minimizes an energy $E_{a,c_o=0,b=0,\alpha,\beta}$ (see Figure \ref{minimalexamples}), while for $c_o>0$, Patnaik \cite{P94} produced an example of a non axially symmetric constant mean curvature annulus having the same boundary as the surfaces $\mathcal{N}_{i}$. For suitable choice of the parameters these surfaces minimize an energy $E_{a,c_o>0,b=0,\alpha,\beta}$.

We next discuss the stability of the nodoidal domains ${\mathcal N}_i$, $i=1,...,4$, when they are critical but not minimizing. 

\begin{prop} For $b>0$ the nodoidal domains ${\mathcal N}_2$ and $\mathcal{N}_3$ are unstable, while for $b<0$, the nodoidal domains ${\mathcal N}_1$, $\mathcal{N}_3$ and $\mathcal{N}_4$ are unstable.
\end{prop}
{\it Proof.\:} We will prove that the nodoidal domain $\mathcal{N}_2$ is unstable for $b>0$, the proof of the other statements is similar.

We fix the value of the constant mean curvature $H=-c_o$ and regard the formulas \eqref{nodc} and \eqref{nodc2} as defining a one parameter family of nodoids with flux parameter $\varpi=\varpi_o+\epsilon$, where $\varpi_o=-c_o\alpha/\beta$ is the flux parameter for ${\mathcal N}_2$. For $\epsilon\geq 0$, we let $\Sigma_\epsilon$ denote the convex annular regions in these nodoids bounded by two circles of radii $r_o:=\sqrt{\alpha/\beta}$ and we let $u^*(\epsilon)$ denote the value of $u$ on $\partial\Sigma_\epsilon$. From $-c_o r_o^2=\varpi_o$ and $r_ou^*(\epsilon)-c_or_o^2=\varpi_o+\epsilon$, we get that $u^*(\epsilon)=\epsilon/r_o$.

The total curvature of $\Sigma_\epsilon$ is easily computed to be
$$\int_{\Sigma_\epsilon}K\:d\Sigma=4\pi\sqrt{1-\left(\frac{\epsilon}{r_o}\right)^2}\:.$$
Therefore, the energy of the domains $\Sigma_\epsilon$ is given by
$$E[\Sigma_\epsilon]=4\pi b\sqrt{1-\left(\frac{\epsilon}{r_o}\right)^2}+8\pi\sqrt{\alpha\,\beta}\,,$$
since $H\equiv -c_o$ and $\alpha\kappa^2=\beta$ holds on the boundary.

Finally, using this expression of the energy in terms of $\epsilon\geq 0$, we conclude that $\partial^2_{\epsilon \epsilon}E[\Sigma_\epsilon]\lvert_{\epsilon=0}=-4\pi b/r_o^2<0$ holds since $b>0$, i.e. $\mathcal{N}_2$ is unstable. {\bf q.e.d.}
\\

\section*{Appendix A. First Variation Formula}

In this appendix, we will compute the first variation formula of the total energy $E$. Most of the calculations below are well known and they are included for completeness.

Let $\delta X$ be a smooth ${\bf R}^3$ valued map on $\Sigma$ which we consider as a variation field, i.e. the linear term of a deformation $X_\epsilon:=X+\epsilon\, \delta X+\mathcal{O}(\epsilon^2)$.

Although for functionals with geometric character (i.e. invariant under changes of parameterization) only normal variations are usually considered, for surfaces with boundary it is essential to use variations having normal as well as tangential components. In general, the boundary is not invariant under tangential variations and, hence, computing the first variation formula  only for normal variations leads to a loss of valuable information.  

We decompose $\delta X$ as $V+\psi\nu$ with $V$ tangent to the surface. At times we will denote with ``dot'' derivatives with respect to the variation parameter $\epsilon$, i.e. $\delta f={\dot f}$ for an $\epsilon$ dependent function $f$ on $\Sigma$.

We first consider the case $V\equiv 0$. A straightforward calculation gives the variation of the normal field as
\begin{equation*}
{\dot \nu}=-\nabla \psi\:,
\end{equation*}   
where we denote the surface gradient operator by $\nabla$. 

We next compute the variation of the components metric tensor $g_{ij}$ with respect to a locally defined frame field $\{e_1, e_2\}$ which is orthonormal for the metric induced by $X_0\equiv X$, i.e. $g_{ij}(0)=\delta_{ij}$,
\begin{eqnarray*}
{\dot g}_{ij}&=& \delta \left(X_i\cdot X_j\right)={\dot X}_i\cdot X_j+X_i\cdot {\dot X}_j=\left(\psi\nu\right)_i\cdot X_j+X_i\cdot \left(\psi\nu\right)_j\\
&=&2\psi \nu_i\cdot X_j=-2\psi L_{ij}\,,
\end{eqnarray*}
where $L_{ij}:=-X_i\cdot \nu_j$ are the components of the second fundamental form of the surface.  Recall that this tensor is symmetric, i.e. $L_{ij}=L_{ji}$.

The induced surface measure on $\Sigma$ is given by $d\Sigma=g^{1/2}\left(e_1^*\wedge e_2^*\right)$, where $g:=\det (g_{ij})$ and $\{e_1^*,e_2^*\}$ is the dual basis. From this and the calculation above, one easily obtains
$$\delta d\Sigma=-2H\psi d\Sigma\,.$$

We now consider the variation of the second fundamental form of the surface $\left(L_{ij}\right)$:
\begin{eqnarray*}
{\dot L}_{ij}&=&- \delta \left(X_i\cdot \nu_j\right)=- {\dot X}_i\cdot \nu_j-X_i\cdot {\dot \nu}_j=-\left(\psi \nu\right)_i\cdot \nu_j-X_i\cdot \left(-\nabla \psi\right)_j\\
&=&-\psi\nu_i\cdot \nu_j+\psi_{,ij}=-\psi L_{ik}L_{kj}+\psi_{,ij}\,.
\end{eqnarray*}

From the previous two calculations, we obtain, using ${\dot g}^{ij}=-{\dot g}_{ij}$,
\begin{eqnarray*}
{\dot H}&=&\delta \left(\frac{1}{2}g^{ij}L_{ij}\right)= \frac{1}{2}\left({\dot g}^{ij}L_{ij}+g^{ij}{\dot L}_{ij}\right)\\&=& \frac{1}{2}\left(2\psi L_{ij}L_{ij}+g^{ij}\psi_{,ij}-\psi L_{ij}L_{ij}\right)=\frac{1}{2}\left(\Delta \psi+\lVert d\nu\rVert^2\psi\right),
\end{eqnarray*}
where $\lVert d\nu\rVert^2=L_{ij}L_{ij}=4H^2-2K$ is the square of the norm of the second fundamental form. 

In order to compute the pointwise variation of the Gaussian curvature $K$, we choose the frame so that, in addition to being orthonormal, $L_{ij}=k_i\delta_{ij}$ holds. Here, $k_1$ and $k_2$ are the principal curvatures of the surface. Then we obtain
\begin{eqnarray*}
{\dot K}&=&\delta \left( \frac{L_{11}L_{22}-L_{12}^2}{g}\right)={\dot L}_{11}k_2+k_1{\dot L}_{22}-{\dot g}K\\
&=&k_2\psi_{,11}+k_1\psi_{,22}-\left(k_2k_1^2+k_1k_2^2\right)\psi+4HK\psi\\
&=&k_2\psi_{,11}+k_1\psi_{,22}+2HK\psi\,.
\end{eqnarray*}
The Codazzi equations with respect to this frame are:
\begin{eqnarray*}
e_2(k_1)=:k_{1,2}&=&-\left(k_1-k_2\right)\nabla_1e_2\cdot e_1\:,\\ 
e_1(k_2)=:k_{2,1}&=&\left(k_1-k_2\right)\nabla _2e_1\cdot e_2\:.
\end{eqnarray*}

Define $A:=\left(d\nu+2H\,Id\right)\nabla \psi= k_2\psi_{1}e_1+k_1\psi_{2}e_2$. Below we use the notation $\psi_{i,j}$ to denote $e_j(\psi_i)$ and as before 
$\psi_{,ij}$ denotes the $(i,j)$ component of the Hessian tensor.  Using the  Codazzi equations, we have that the divergence $\nabla\cdot$ of $A$ is given by
\begin{eqnarray*}
\nabla\cdot A &=& \nabla_i\left( k_2\psi_{1}e_1+k_1\psi_{2}e_2\right)\cdot e_i\\
&=&\left(k_2\psi_1\right)_1+\left(k_1\psi_2\right)_2+k_2\psi_1\nabla_2e_1\cdot e_2+k_1\psi_2\nabla_1e_2\cdot e_1\\
&=&k_{2,1}\psi_1+k_2\psi_{1,1}+k_{1,2}\psi_2+k_1\psi_{2,2}+k_2\psi_1\nabla_2e_1\cdot e_2+k_1\psi_2\nabla_1e_2\cdot e_1\\
&=&\psi_1\left(k_1-k_2\right)\nabla _2e_1\cdot e_2+k_2\psi_{1,1}-\psi_2\left(k_1-k_2\right)\nabla_1e_2\cdot e_1\\
&&+k_1\psi_{2,2}+k_2\psi_1\nabla_2e_1\cdot e_2+k_1\psi_2\nabla_1e_2\cdot e_1\\
&=&\psi_1k_1\nabla _2e_1\cdot e_2+k_2\psi_{1,1}+\psi_2k_2\nabla_1e_2\cdot e_1+k_1\psi_{2,2}\\
&=&k_1\psi_{,22}+k_2\psi_{,11}\:.
\end{eqnarray*}
Comparing this with the expression for ${\dot K}$ above, we obtain
$${\dot K}=\nabla\cdot \left(\left[d\nu+2H\, Id\right]\nabla \psi\right)+2HK\psi\:.$$

For the case of a variation tangential to the surface, i.e. $\delta X=V$, it is clear that ${\dot H}=\nabla H\cdot V$, $\dot K=\nabla K\cdot V$ and $\delta d\Sigma =\left(\nabla \cdot V\right)d\Sigma$, where $\nabla\cdot V$ denotes the divergence of $V$.

We are now in a position to compute the variation of the Helfrich energy (\cite{H})
$$\mathcal{H}[\Sigma]=\int_\Sigma a\left(H+c_o\right)^2\,d\Sigma\,.$$ 
For $\delta X=\psi\nu +V$, integrating by parts we obtain
\begin{eqnarray*}
\delta \mathcal{H}[\Sigma]&=& \int_\Sigma 2\left(H+c_o\right){\dot H}\:d\Sigma +\int_\Sigma \left(H+c_o\right)^2\delta(d\Sigma)\\
&=& \int_\Sigma 2\left(H+c_o\right)\left(\frac{1}{2}\left[\Delta \psi+(4H^2-2K)^2\psi\right] +\nabla H\cdot V \right)d\Sigma\\
&&+\int_\Sigma \left(H+c_o\right)^2 \left(-2H\psi+\nabla\cdot V\right)d\Sigma\\
&=& \int_\Sigma \left(H+c_o\right)\left( \Delta \psi+2(H^2-K-c_oH)\psi \right)+\nabla \cdot \left[(H+c_o)^2V\right]d\Sigma\\
&=&\int_\Sigma \left[ \Delta(H+c_o)+2(H+c_o)(H^2-K-c_oH) \right]\psi\,d\Sigma\\
&&+\oint_{\partial \Sigma} \left(H+c_o\right)\partial_n\psi- \partial_n\left(H+c_o\right)\psi+\left(H+c_o\right)^2V\cdot n\:ds\:.
\end{eqnarray*}
In the last line, we have used Green's Second Identity.

Similarly, for $\delta X=\psi\nu+V$ the variation of the total Gaussian curvature is given by
\begin{eqnarray*}
\delta\left(\int_\Sigma K\:d\Sigma\right)&=&\int_\Sigma \dot{K}\,d\Sigma+K\,\delta(d\Sigma)\\
&=&\int_\Sigma \left( \nabla\cdot \left(\left[d\nu+2H\, Id\right]\nabla\psi\right)+\nabla K\cdot V+K\nabla\cdot V\right)d\Sigma \\
&=&\oint_{\partial \Sigma} \left(\left[d\nu+2H\, Id\right]\nabla\psi+K V\right)\cdot n\:ds\\
&=&\oint_{\partial \Sigma} \left(\kappa_n\partial_n\psi-\tau_g\psi'+K V\cdot n\right)ds\,,
\end{eqnarray*}
where in the third line we have used the Divergence Theorem.

Finally, to compute the first variation formula of the bending energy of the boundary we need to know the pointwise variation of the squared curvature. For an arbitrary parameter $t$, we have
\begin{eqnarray*}
\delta\left(\kappa^2\right)&=&\delta\left(T'\cdot T'\right)=\delta\left(\frac{T_t}{\lVert C_t\rVert}\cdot\frac{T_t}{\lVert C_t\rVert}\right)=2\delta\left(\frac{T_t}{\lVert C_t\rVert}\right)\cdot T'\\
&=&2\left(\dot{T}\right)'\cdot T'-2\kappa^2 T\cdot \left(\dot{C}\right)'=2T' \cdot \left(\dot{C}\right)''-4\kappa^2 T\cdot \left(\dot{C}\right)',
\end{eqnarray*}
since $\dot{T}=\left[\left(\dot{C}\right)'\right]^\perp$ and, hence, $\left(\dot{T}\right)'=\left(\dot{C}\right)''-T'\cdot\left(\dot{C}\right)'-T\cdot\left(\dot{C}\right)''$. Next, from the induced measure on $\partial\Sigma$, we directly obtain $\delta(ds)=T\cdot\left(\dot{C}\right)'ds$. Combining both things we get
\begin{eqnarray*}
\delta \left( \oint_{\partial \Sigma}\left[ \alpha \kappa^2+\beta\right]ds\right)&=&\oint_{\partial\Sigma} \alpha\delta\left(\kappa^2\right)ds+\left(\alpha\kappa^2+\beta\right)\delta(ds)\\
&=&\oint_{\partial\Sigma} \left(2\alpha T'\cdot \left[\dot{C}\right]''-\left[3\alpha\kappa^2-\beta\right]T\cdot\left[\dot{C}\right]'\right)ds\\
&=&\oint_{\partial\Sigma}\left(2\alpha T'''+\left[(3\alpha\kappa^2-\beta)T\right]'\right)\cdot \dot{C}\,ds\,,
\end{eqnarray*}
where in the last equality we have integrated by parts.

The same variations have been obtained using different techniques in \cite{BMF} and \cite{OT1}, to mention a couple. For the boundary energy see also the Appendix of \cite{LS}.

\section*{Acknowledgments}

The second author has been partially supported by MINECO-FEDER grant PGC2018-098409-B-100, Gobierno Vasco grant IT1094-16 and by Programa Posdoctoral del Gobierno Vasco, 2018. He would also like to thank the Department of Mathematics and Statistics of Idaho State University for its warm hospitality.

\bigskip

\begin{flushleft}
Bennett P{\footnotesize ALMER}\\
Department of Mathematics,
Idaho State University,
Pocatello, ID 83209,
U.S.A.\\
E-mail: palmbenn@isu.edu
\end{flushleft}

\bigskip

\begin{flushleft}
\'Alvaro P{\footnotesize \'AMPANO}\\
Department of Mathematics, University of the Basque Country, Bilbao, Spain.\\
Department of Mathematics, Idaho State University,
Pocatello, ID 83209,
U.S.A.\\
E-mail: alvaro.pampano@ehu.es
\end{flushleft} 

\end{document}